\documentclass[amsart]{article}

\usepackage{amsmath}
\usepackage{amssymb}
\usepackage{amscd}
\usepackage{bbm}
\usepackage{cyrillic}
\usepackage{latexsym}
\usepackage{makeidx}
\usepackage{mathrsfs}
\usepackage{mathtools}
\usepackage{theorem}
\usepackage[all]{xy}
\usepackage[colorlinks]{hyperref}
\hypersetup{colorlinks=false , linkcolor=blue}

\DeclareMathAlphabet{\mathbbold}{U}{bbold}{m}{n}

\newtheorem{teo}{Theorem}[section]
\newtheorem{propo}[teo]{Proposition}
\newtheorem{lema}[teo]{Lemma}
\newtheorem{coro}[teo]{Corollary}

{\theorembodyfont{\rmfamily }
\newtheorem{defi}[teo]{Definition}

\newtheorem{obser}[teo]{Remark}
\newtheorem{ejem}[teo]{Example}
\newtheorem{parra}[teo]{}}
\def\demo{\noindent \textit{Proof: }}
\def\fin{\end{document}}

\DeclareMathOperator{\codim}{codim}

 \DeclareMathOperator{\hofib}{hofib}
\DeclareMathOperator{\Hom}{Hom} 
 \DeclareMathOperator{\Pic}{Pic}
 \DeclareMathOperator{\rk}{rank}
 \DeclareMathOperator{\Spec}{Spec}
\DeclareMathOperator{\Td}{Td} \DeclareMathOperator{\Th}{Th}

\newcommand{\cyrrm}{\fontencoding{OT2}\selectfont\textcyrup}
\def\11{1}
\def\11{\mathbbm{1}}

\def\ad{\mathrm{ad}}

\def\cf{\emph{cf. }}
\def\ch{\mathrm{ch}}
\def\cl{\mathrm{cl}}

\def\d{\mathrm{d}}

\def\Fr{\mathrm{Fr}}

\def\hatch{\widehat{\mathrm{ch}}}
\def\hatHm{\widehat{ \mathrm{H}}_{\cyrrm{B}}}
\def\hatKH{\widehat{KH}}

\def\new{\mathrm{new}}

\def\p{\mathfrak{p}}

\def\raya{\ \underline{\phantom{a}}\ }
\def\syn{\mathrm{syn}}
\def\t{\mathbf{t}}

\def\qed{\hspace*{\fill }$\square $ }

                    \def\AA{\mathbb{A}}   

                    \def\CC{\mathbb{C}} 
                     
 \def\EE{\mathbb{E}}
\def\F{\mathcal{F}} \def\FF{\mathbb{F}} 
 \def\GG{\mathbb{G}}
                                         \def\HH{\mathbf{H}}
\def\I{\mathcal{I}}
\def\K{\mathcal{K}}                     
\def\L{\mathcal{L}}                     
 \def\NN{\mathbb{N}}
\def\M{\mathcal{M}} \def\MM{\mathbb{M}}
\def\O{\mathcal{O}}
 \def\PP{\mathbb{P}}
                    \def\QQ{\mathbb{Q}}
                    \def\RR{\mathbb{R}}

\def\V{\mathcal{V}}

 \def\ZZ{\mathbb{Z}}

\def\MGL{\mathrm{MGL}}
\def\KGL{\mathrm{KGL}}

\def\HB{\mathrm{H}_{\cyrrm{B}}}
\def\hatHB{\widehat{\mathrm{H}}_{\cyrrm{B}}}
\def\hatKGL{\widehat{\mathrm{KGL}}}

\def\DMb{\mathbf{DM}_{\cyrrm{B}}}
\def\Ho{\mathbf{Ho}}
\def\mod{\mbox{-}\mathbf{mod}}
\def\Mon{\mathbf{Mon}}

\def\SH{\mathbf{SH}}

\def\Sm{\mathbf{Sm}}

\def\Spt{\mathbf{Spt}}

\begin{document}

\title{Riemann-Roch for homotopy invariant $K$-theory and Gysin morphisms}

\author{A. Navarro \thanks{ Instituto de Ciencias Matem\'{a}ticas (CSIC-UAM-UC3M-UCM) }}

\date{May 3$^{\mathrm{th}}$, 2016}

\maketitle

\begin{abstract}
We prove the Riemann-Roch theorem for homotopy invariant $K$-theory and projective
local complete intersection morphisms between finite dimensional noetherian schemes,
without smoothness assumptions. We also prove a new Riemann-Roch theorem for the
relative cohomology of a morphism.

In order to do so, we construct and characterize Gysin morphisms for regular
immersions between cohomologies represented by spectra (examples include homotopy
invariant $K$-theory, motivic cohomology, their arithmetic counterparts, real
absolute Hodge and Deligne-Beilinson cohomology, rigid syntomic cohomology, mixed
Weil cohomologies) and use this construction to prove a motivic version of the
Riemann-Roch.
\end{abstract}

\tableofcontents

\newpage

\section*{Introduction}

The original Grothendieck's Riemann-Roch theorem states that for any proper morphism
$f\colon Y\to X$, between nonsingular quasiprojective irreducible varieties over a
field, and any element $a \in K_0(Y)$ of the Grothendieck group of vector bundles the
relation
$$
\ch(f_!(a))= f_*\bigl(\Td(T_f)\cdot \ch(a)\bigr)
$$
holds (\emph{cf.} \cite{BorelSerre}). Recall that $\ch$ denotes the Chern character,
$\Td (T_f)$ the Todd class of the relative tangent bundle and $f_*$ and $f_!$ the
direct image in the Chow ring and $K_0$ respectively. Later the result was
generalized for locally complete intersection morphisms between singular projective
algebraic schemes (\cf \cite{SGA6}, \cite{BFM}).

The extension to higher $K$-theory and schemes over a regular base was proved by
Gillet in \cite{Gillet}. The Riemann-Roch theorem proved there is for projective
morphisms between smooth quasiprojective schemes. Therefore, note that in the case
over a field Gillet's theorem does not recover the result of \cite{SGA6}. The
furthest generalization of the higher Riemann-Roch theorem I know is \cite{Deglise}
and \cite{HS} where D\'{e}glise and Holmstrom-Scholbach independently obtained the
Riemann-Roch theorem for higher $K$-theory and projective morphisms between regular
schemes.

After the work of Cisinski in \cite{Cisinski} on Weibel's homotopy invariant
$K$-theory one may apply Voevodsky's motivic homotopy theory to it. With Cisinski's
result, we present a Riemann-Roch theorem for homotopy invariant $K$-theory $KH$ and
projective lci morphisms without smoothness assumptions on the schemes. More
concretely, the theorem we prove for motivic cohomology $H_{\M}$ and finite
dimensional noetherian schemes is the following:

\medskip

\noindent \textbf{Theorem:} \emph{ Let $f\colon Y\to X$ be a projective lci morphism
and denote $T_f\in K_0(Y)$ the virtual tangent bundle and $\Td$ the multiplicative
extension of the series given by $\frac{t}{1-e^{-t}}$. Then the diagram
$$
\xymatrix{
KH(Y)_\QQ\ar[r]^{f_*} \ar[d]_{\Td(T_f)\ch}& KH(X)_\QQ\ar[d]^{\ch}\\
H_{\M}(Y, \QQ) \ar[r]^{ f_*} & H_{\M} (X, \QQ) }
$$
commutes. In other words, for $a \in KH(Y)_\QQ$ we have
$$
\ch(f_*(a))= f_*\bigl(\Td(T_f)\cdot \ch(a)\bigr).
$$}

\medskip

\noindent From here, we deduce Riemann-Roch theorems for many cohomologies. In
particular, for real absolute Hodge and Deligne Beilinson cohomology, rigid syntomic
cohomology and mixed Weil cohomologies such as algebraic de Rham and geometric \'{e}tale
cohomology. Note also that, due to the comparison result with \emph{cdh}-motivic
cohomology in \cite{CD3}, motivic cohomology on singular schemes over a base field of
exponential characteristic is computed by explicit cycles.

In order to prove this result, the addition we make to the theory is the construction
of the Gysin morphism for regular immersions and every cohomology given by spectra.
Since its beginnings, the standard construction of the Gysin morphism in motivic
homotopy theory relies on the Thom space and the purity isomorphism. However, purity
requires smoothness assumptions (\emph{cf.} \cite{MV}). Our approach is a different
one: We lift the work of Gabber for \'{e}tale cohomology to the motivic homotopy setting
and, thus, obtain Gysin morphisms for regular immersions without smoothness
assumptions on the schemes. This leads to the construction of new Gysin morphisms for
many theories like homotopy invariant $K$-theory, motivic cohomology, real absolute
Hodge and Deligne-Beilinson cohomology, rigid syntomic cohomology, motivic
\emph{cdh}-cohomology and any cohomology coming from a mixed Weil theory.

Modules over a cohomology theory are notable geometric and arithmetic invariants.
Recall that higher arithmetic $K$-theory and arithmetic motivic cohomology are
modules over $K$-theory and motivic cohomology respectively. In addition, for any
cohomology the relative cohomology of a morphism is also a module. Note that the
cohomology with proper support, the cohomology with support on a closed subscheme,
and the reduced cohomology are the relative cohomology of a closed immersion, an open
immersion and the projection over a base point respectively.

We deduce from our results new Gysin morphisms and a Riemann-Roch theorem for
abstract modules. In particular, we prove a new Riemann-Roch theorem for relative
cohomology and finite dimensional noetherian schemes.

\medskip

\noindent \textbf{Theorem:} \emph{Let $g\colon T \to X$ be a morphism of schemes,
$f\colon Y\to X$ be a projective lci morphism. Denote $g_{{\scriptscriptstyle
Y}}\colon T\times_X Y\to Y$, $T_f\in K_0(Y)$ the virtual tangent bundle of $f$ and
$KH(g)$ and $H_{\M}(g,\QQ)$ the relative homotopy invariant $K$-theory and motivic
cohomology of $g$ respectively. Assume in addition either $g$ is proper or $f$ is
smooth, then the diagram
$$
\xymatrix{
KH(g_{{\scriptscriptstyle Y}})_\QQ\ar[r]^{f_*} \ar[d]_{\Td(T_f)\ch}& KH(g)_\QQ\ar[d]
^{\ch}\\
H_{\M}(g_{{\scriptscriptstyle Y}},\QQ) \ar[r]^{f_*} & H_{\M}(g,\QQ) }
$$
commutes. In other words, for $m \in KH(g_{\scriptscriptstyle Y})_\QQ$ we have}
$$
\ch(f_*(m))=  f_*\bigl(\Td(T_f)\cdot \ch(m)\bigr).
$$

\medskip

\noindent We also deduce the higher arithmetic Riemann-Roch theorem of \cite{HS} from
the Riemann-Roch theorem for modules.

\medskip

The paper is organized as follows: Section 1 recalls basic facts of the stable
homotopy category and cohomologies defined by spectra, in particular, orientations or
Chern classes. In section \ref{Sec Ejemplos Modulos} we introduce the relative
cohomology in the context of motivic stable homotopy theory and construct an absolute
spectrum which represents it under some conditions. We also recall
Holmstrom-Scholbach's construction of the higher arithmetic $K$-theory and motivic
cohomology. In section 2 we construct Gysin morphisms using Gabber's ideas for the
case of regular immersions and prove the basic properties. We prove the motivic
Riemann-Roch theorem in section 3 and an analogous statement for modules (\emph{cf.}
Theorem \ref{Motivic Riemann Roch} and Theorem \ref{Teo RR modulos}) and obtain the
Riemann-Roch theorem for homotopy invariant $K$-theory and for relative cohomologies
as a result. The Appendix is devoted to the explicit construction of the real
absolute Hodge spectrum using Burgos's complex (\emph{cf.} \cite{Burgos3}) and to
check that this spectrum, as well as the Deligne-Beilinson spectrum of \cite{HS},
represents their cohomologies also in the singular context.

\subsection*{Acknowledgements}
I am deeply thankful to my advisors, J.I. Burgos Gil and F. D\'{e}glise, for guiding me
and generously sharing with me their knowledge during the preparation of this work. I
also want to thank the geometers of Departamento de Matem\'{a}ticas at UEX, the time I
spent with them and our discussions greatly benefited this work. I would also like to
thank J. Riou for many fruitful remarks on a previous version of this work.

The author has been partially supported by ICMAT Severo Ochoa project SEV-2015-0554
(MINECO), MTM2013-42135-P (MINECO) and ANR-12-BS01-0002 (ANR).

\section{Stable homotopy}

All schemes considered are noetherian, finite dimensional and over a fixed base
scheme $S$.

\medskip

Let $X$ be a scheme and denote $\SH(X)$ the stable homotopy category over $X$
(\cite[$\S$ 5]{Voevodsky}). The categories obtained if one varies $X$ satisfies
Grothendieck's six functor formalism. We will only recall here some of their
properties that will be used later. We refer to \cite{Ayoub} and \cite{CD} for a
complete exposition of the construction of $\SH$ and Grothendieck's six operations in
this context.

Let $\Sm _X$ be the category of smooth schemes over $X$. Recall that the category
$\SH(X)$, whose objects are called \emph{spectra}, is constructed out of the category
of pointed Nisnevich presheaves of simplicial sets on $\Sm_X$. Any $X$-scheme $T$
defines the simplicial sheaf given by $\Hom_{\Sm _X}$ $(\raya  , T)$, which we still
denote $T$. We denote $T_+=T\sqcup \mathit{pt}$ the pointed simplicial sheaf it
defines. Then we have the infinite suspension functor
$$
\Sm _X \xrightarrow{\Sigma^\infty} \SH (X).
$$
We denote $\Sigma^\infty T$, or simply $T$ if no confusion is possible. In the
particular case of $X$ we will use the notation $\11_X\in \SH (X)$.

The category $\SH (X)$ is a monoidal triangulated category. We denote the shift
functor as
$$
[1]\colon\SH (X)\to  \SH (X)
$$
and the product is denoted by $\wedge $. Recall that the wedge product is defined in
the category of \emph{symmetric} spectra $\Spt_\Sigma(X)$ whose localization by a
suitable model structure is equivalent to $\SH(X)$.

Recall that we consider the projective line pointed at the infinity. The invertible
element $\11_X(1) \coloneqq \mathrm{coker} (\Sigma^\infty X\to
\Sigma^{\infty}\PP^1)[-2]$ is called the \emph{Tate object}. For any spectrum $E$ we
denote the \emph{Tate twist} by $E(1)\coloneqq E\wedge \11_X (1)$ and denote
$E(q)[p]$ for twisting and shifting $q$, $p \in \ZZ$ times respectively.

For any morphism of schemes $f\colon Y \to X$ there is a pullback functor
$$
f^*\colon \SH (X) \to \SH (Y)
$$
which is also functorial with respect to $f$. Moreover, the functor $f^*$ admits a
right adjoint denoted $f_*$.

The pullback functor $f^*$ is also monoidal and therefore $f^*(\11_Y)=\11_X$. If
$i\colon Y\to X$ is a closed immersion and $U=X-Y$, then $\11_Y\simeq
i^*(\Sigma^\infty X/U)$ where $\Sigma^\infty X/U$ stands for the spectrum given by
the pointed quotient Nisnevich presheaf $X/U$. Moreover, the pullback and pushforward
functor satisfy the projection formula. In other words, for every $E\in \SH(Y)$ and
$F\in \SH(X)$ there is a canonical isomorphism
$$
i_*(E\wedge i^* (F))\xrightarrow\sim i_*(E)\wedge F.
$$

The stable homotopy category relates to other categories that we will use at some
point. In \cite{MV} Morel and Voevodsky defined the \emph{homotopy category} of
schemes $\HH (X)$. This category is constructed as a localization of an intermediate
category $\HH_s (X)$ called the \emph{simplicial homotopy category} of schemes. Their
pointed counterparts are related to the stable homotopy category via the derived
infinite suspension so that there are functors
$$
\HH_\bullet^s (X)\to \HH_\bullet (X)
\buildrel{\Sigma^\infty}\over{\longrightarrow}\SH (X).
$$
As an analogy with the topological case, the $\Hom$ in the homotopy category
$\HH_\bullet$ is denoted $[\raya , \raya ]$.

Finally, recall that the category of Beilinson motives $\DMb (X)$ constructed in
\cite[$\S$ 14]{CD} satisfies Grothendieck's six functor formalism and it is a full
subcategory of $\SH(X)_\QQ$.

\subsection{Cohomology and its operations.}

\begin{defi}\label{defi absoluto}
A \textbf{ring spectrum} is an associative commutative unitary monoid in the stable
homotopy category $\SH (X)$ of a scheme $X$. Thus, a ring spectrum is a triple $(E
,\, \mu\colon E\wedge E\to E ,\, \eta\colon X \to E)$ consisting of a spectrum, the
product morphism, and the unit morphism, satisfying the usual conditions of
commutative monoids. A \textbf{morphism} of ring spectra is morphism of spectra
compatible with the unit and product morphism.

An \textbf{absolute ring spectrum} $\EE$ is a family of ring spectra $\EE_X$ on
$\SH(X)$ for every scheme $X$ stable by pullback. In other words, for every morphism
$f\colon Y\to X$ we have fixed an isomorphism $\epsilon_f\colon f^*\EE_X \to \EE_Y$
satisfying the usual compatibility conditions. A \textbf{morphism} of absolute ring
spectra $\varphi \colon \EE \to \FF$ is a collection of morphism of ring spectra
$\varphi_X \colon \EE_X \to \FF_X$ for every scheme $X$ stable by pullback. In other
words, for every morphism $f\colon Y \to X$ we have $ \epsilon _f^\FF\circ
f^*(\varphi_X)= \varphi_Y\circ \epsilon_f^\EE.$ Let $X$ be a scheme. We call a family
stable by pullback of spectra $ \EE_Y $ for $Y$ an $X$-scheme an absolute spectra
\emph{over $X$}.

Let $E$ in $\SH(X)$ be a ring spectrum, an \textbf{$E$-module} is a spectrum $M$ in
$\SH(X)$ together with a morphism of spectra $\upsilon\colon E\wedge M\to M$ in
$\SH(X)$ satisfying the usual module condition. Let $\varphi\colon E\to F$ be a
morphism of ring spectra and $(M' , \, v')$ be an $F$-module. A
\textbf{$\varphi$-morphism} of modules $\Phi \colon M\to M'$ is a morphism of spectra
in $\SH(X)$ such that $ v'\wedge \varphi \wedge \Phi=\Phi\circ v$.

Let $\EE$ be an absolute ring spectrum, an \textbf{absolute $\EE$-module} is an
absolute spectrum $\MM$ such that $\MM_X$ is a $\EE_X$-module for every $X$ and the
structural isomorphisms $\epsilon_f$ are isomorphisms of modules. Let $\varphi\colon
\EE \to \FF$ be a morphism of absolute ring spectra and $\MM'$ be an absolute
$\FF$-module. A $\varphi$-\textbf{morphism} of absolute modules $\Phi\colon \MM\to
\MM'$ is a morphism of absolute spectra such that $\Phi_X$ is a $\varphi_X$-morphism
of modules which are stable by pullback. On the following, we may omit the adjective
\emph{absolute} when it is clear by the context and notation.
\end{defi}

\begin{obser}
Every absolute spectrum $\EE$ is naturally isomorphic to the absolute spectrum
obtained by pullback of the spectrum $\EE_S\in \SH(S)$, where $S$ is the base scheme.
Instead of considering schemes over a fixed base $S$, one may work over a general
category of schemes $\mathbf{S} $ without a final object. All definitions and proofs
of this paper may be carried into this context using the notion of
$\mathbf{S}$-absolute spectra (\emph{cf.} \cite{Deglise}). However, since we will not
make use of this generality we have chosen otherwise. In addition, we will abuse
notation and call $\EE_S$ also the absolute spectrum.
\end{obser}

\begin{defi}
Let $X$ be a scheme and $\EE$ be an absolute spectrum over $X$, we define the
\textbf{$\EE$-cohomology} of $X$ to be
$$
\EE ^{p,q} (X)=\Hom _{\SH(X)}(\11_X,\EE_X (q)[p])\qquad \mbox{for }(p,q)\in \ZZ^2
$$
and $\EE (X)=\bigoplus_{p,q} \EE ^{p,q}(X)$. Let $i\colon Z\to X$ be a closed
immersion, we define the \textbf{$\EE$-cohomology with support on $Z$} to be
$$
\EE ^{p,q}_Z (X)=\Hom _{\SH(X)}(i_*\11_Z,\EE_X(q)[p]) \quad \mbox{for }p,q\in \ZZ.
$$

For any $f\colon  T \to X$ and any closed immersion $i\colon Z\to X$ we define the
\textbf{inverse image of }$f$ which maps an element $a\colon i_*\11_Z\to \EE_X(q)[p]$
in $ \EE_Z^{p,q} (X)$ to the composition
$$
f^*(a)\colon f^*i_*\11_Z \simeq i'_*\11_{Z'} \rightarrow f^*(\EE_{X}(q)[p])\simeq
\EE_{T} (q)[p]\ \in \EE_{Z'}^{p,q}(T)
$$
where $i'\colon Z'=Z\times_{T}  X\to T$. We denote it $f^*\colon \EE_Z(X)\to
\EE_{Z'}(T)$.
\end{defi}

\begin{ejem}\label{Ejem espectros en anillo}
\begin{itemize}

\item
Let $k$ be a perfect field and consider $S=\Spec (k)$. In \cite[2.1.4]{CD2} Cisinski
and D\'{e}glise defined the notion of mixed Weil theories with coefficient in a field of
characteristic zero. In \cite{CD2}, every such theory is proved to define a ring
spectrum on $\SH(S)_\QQ$ stable by pullback and, therefore, an absolute ring
spectrum. Recall that there is an algebraic de Rham ring spectrum for $k$ of
characteristic zero, an analytic de Rham ring spectrum for $k$ algebraically closed
of characteristic zero, and a $\QQ_l$ geometric \'{e}tale ring spectrum for $k$ countable
perfect and $l$ a prime different from the characteristic of $k$.

\item
Let $S=\Spec (\ZZ)$. The $K$-theory absolute ring spectrum $\KGL$ is defined in
\cite{Voevodsky} and \cite{Riou3}. By constructions it is periodic, meaning that
there are isomorphisms
$$
\KGL \simeq \KGL (i)[2i] \, , \, \forall i.
$$
It represents Weibel's homotopy invariant $K$-theory for every scheme (\cf
\cite[2.15]{Cisinski}), and therefore represents Quillen's algebraic $K$-theory for
regular schemes. We denote the cohomology groups they define as $KH _i (\raya )$.

Following \cite[5.3]{Riou3}, the $\QQ$-localization of the $K$-theory spectrum admits
a decomposition induced by the Adams operations, i.e., $ \KGL_\QQ = \bigoplus_{i\in
\ZZ} \KGL_ \QQ ^{(n)}, $ where $\KGL_ \QQ ^{(n)}$ denotes the eigenspaces for the
Adams operations. The Beilinson's absolute motivic cohomology spectrum is defined as
$\HB=\KGL_\QQ^{(0)} \in \SH(S)_\QQ$ and it is also an absolute ring spectrum.
Therefore the Adams operations define an isomorphism
$$
\ch \colon \KGL_\QQ \xrightarrow \sim \bigoplus_{i\in \ZZ}\HB(i)[2i]$$ (see
\cite[5.3.17]{Riou3}). We call this morphism the \textbf{Chern character} since for
any regular scheme $X$ it induces the classical higher Chern characters
$$
\ch_{r,n}\colon  K_r(X)_ \QQ \to H_{\mathcal{M}}^{2n-r}(X,\QQ(n)).
$$
The Chern character is a morphism of absolute ring spectra. In particular $\HB$ is a
module over $\KGL_\QQ$.

Finally, Voevodsky's absolute algebraic cobordism ring spectrum $\MGL$ is constructed
out of the Thom spaces of the canonical bundle of Grassmannians (\cite{Voevodsky})
and its cohomology is called algebraic cobordism.

\item
In \cite{Spitz} Spitzweck defines the absolute motivic cohomology ring spectrum
$\mathrm{H}_{\Lambda}$ with coefficients in $\Lambda$ for schemes over a Dedekind
domain $S$. Over a field, this spectrum coincides with the motivic Eilenberg-MacLane
spectrum, so it represents motivic cohomology. Rationally Spitzweck's spectrum
coincides with the Beilinson's motivic cohomology spectrum $\HB$. For coherence with
notations in \cite{BMS} we denote the motivic cohomology groups as
$H_{\mathcal{M}}^{*}(\raya , \Lambda (*))$ for motivic cohomology with coefficients
in $\Lambda$.

\item
In \cite{CD3} Cisinski and D\'{e}glise defined an absolute ring spectrum representing the
\emph{cdh}-motivic cohomology in $\SH (\Spec k)_{\ZZ[1/p]}$ for $k$ a field of
exponential characteristic $p$. This spectrum is isomorphic to the motivic cohomology
spectrum.

\item
Let $S=\Spec (k)$ for $k$ an arithmetic field (\emph{cf}. Appendix \ref{Defi cuerpo
aritmetico}). In his thesis, Drew constructed the absolute ring spectrum representing
absolute Hodge cohomology with rational coefficients. His construction also holds for
any subfield of the real numbers \cite[2.1.8]{Drew}. In \cite{HS} Holmstrom and
Scholbach defined the real Deligne-Beilinson ring spectrum $\EE_{\mathrm{D},X}\in
\SH(X)_{\QQ}$ for $X$ a smooth $S$-scheme. The absolute ring spectrum
$\EE_{\mathrm{D}, \, S}$ represents the Deligne-Beilinson cohomology with real
coefficients. We include in Appendix \ref{Hodge absoluto} another direct construction
of the real Hodge $\EE_{\mathrm{AH}}$ and the real Deligne-Beilinson absolute ring
spectra.

\item
Let $K$ be a $p$-adic field, $k$ its residue field and $S=\Spec (k)$. The absolute
rigid syntomic ring spectrum $\EE_{\syn}\in \SH(S)_\QQ$, which represents Besser's
rigid syntomic cohomology, is constructed in \cite{DM}.

\item
Denote $\EE$ the absolute ring spectrum defined by a mixed Weil theory as in
\cite{CD2} (for example, algebraic and analytic de Rham or $\QQ_l$ geometric \'{e}tale).
Cisinski and D\'{e}glise constructed a morphism of absolute ring spectra
$$
\cl \colon \HB\to \EE.
$$
We call this morphism the \textbf{cycle class}. Therefore any mixed Weil spectrum
$\EE$ is a module over $\HB$. Note that there is an analogous construction of the
cycle class for Besser's rigid syntomic, absolute Hodge and Deligne-Beilinson spectra
(\cf \cite{DM} and \cite{HS}).

\end{itemize}

\end{ejem}

\begin{obser}
In the reference \cite{HS} the Deligne-Beilinson absolute spectrum is proved to
represent the real Deligne-Beilinson cohomology on smooth schemes asking explicitly
for the nonsmooth case. In the Appendix \ref{Hodge absoluto} we check that it also
represent the real Deligne-Beilinson cohomology for general schemes.
\end{obser}

\begin{defi}\label{Defi producto refinado}
Let $\MM$ be an absolute $\EE$-module and denote $\mu\colon \EE\wedge \MM\to \MM $
the structure morphism. Let $Z \xrightarrow j Y \xrightarrow i X$ be closed
immersions. We call the \textbf{refined product} to the morphism
$$
\begin{matrix}
\MM^{p,q}_Z(Y) \times \EE^{r,s}_Y (X) & \to & \MM^{p+r,q+s}_Z (X)\\
(m, a) & \mapsto & m \cdot a
\end{matrix}
$$
defined as follows: let $m \colon j_* (\11_Z) \to \MM_ Y
 (q)[p] $ and $a \colon i_*\11_Y \to
\EE_{X} (s)[r],$ then
$$
m\cdot a\colon i_*j_*\11_Z\longrightarrow i_*\MM_Y(q)[p]\simeq \MM_X\wedge
i_*\11_Y(q)[p]\buildrel{\mathrm{id}\wedge a}\over{\longrightarrow}\MM_X\wedge
\EE_{X}(q+s)[p+r]
$$
$$
\buildrel{\mu}\over{\longrightarrow}\MM_X(q+s)[p+r].
$$
Note that the same construction defines a product $ \EE^{p,q}_Z(Y) \times \MM^{r,s}_Y
(X) \longrightarrow \MM^{p+r,q+s}_Z (X)$.
\end{defi}

Finally, let us recall a generalization of the morphism of forgetting support.

\begin{defi}
Let $\MM$ be an absolute $\EE$-module. Consider $Z\xrightarrow j Y \xrightarrow i X$
closed immersions, we define a morphism
$$
j_\flat  \colon\MM_Z (X) \to \MM_Y(X)
$$
as follows. The adjunction morphism $\ad \colon \11_Y \to j_*j^* \11_Y$ defines a
morphism $i_*(\ad) \colon i_* (\11_Y)\to i_*j_*j^*\11_Y=(ij)_*\11_Z$. Let $a\colon
(ij)_*\11_Z\to \MM_X $ be in $\MM_Z (X)$, we define
$$
j_\flat(a)\coloneqq i_* (\11_Y)\xrightarrow {i_*(\ad)}(ij)_*\11_Z\xrightarrow
{a}\MM_X \, \in \MM_Y(X).
$$
\end{defi}

The properties one may expect from the morphism of forgetting support and the former
product are summarized in the following result, which comes from
\cite[1.2.9]{Deglise} for ring spectra.

\begin{propo}\label{Propo propiedades producto}
Let $\EE$ be an absolute ring spectrum and $\MM$ be an absolute $\EE$-module:

\begin{enumerate}
\item
If $j\colon Z\to Y$ and $i\colon Y\to X$ are closed immersion then $ i_\flat j_\flat
= (ij)_\flat$.

\item
Consider the cartesian squares
$$
\xymatrix{Z'\ar[r]^{j'}\ar[d]& Y'\ar[r] \ar[d]^g & X'\ar[d]^f\\
Z\ar[r]^j& Y \ar[r]& X\\}
$$
where the horizontal arrows are closed immersions. Then, for any $\rho \in \MM_Z(X)$
we have $f^*j_\flat (\rho)= j'_\flat f^* (\rho)$.

\item
With the preceding notations, for any pair $(a,\, m)\in \EE _Z (Y) \times \MM_Y (X)$
we have $f^*(a\cdot m) = g^*(a)\cdot f^*(m)$ and $j_\flat (a \cdot m)=j_\flat
(a)\cdot m$. Analogous formulas hold classes in $\MM _Z (Y) \times \EE_Y (X)$.

\item
Consider closed immersions $T\to Z \to Y \to X$. Then for any triple $(a,\, b,\, m)
\in \EE_T(Z)\times \EE_Z(Y)\times \MM_Y(X)$ we have $a \cdot (b \cdot m)=(a \cdot
b)\cdot m$. Analogous formulas hold for classes in $\EE_T(Z)\times \MM_Z(Y)\times
\EE_Y(X)$ and $\MM_T(Z)\times \EE_Z(Y)\times \EE_Y(X)$.

\item
Consider the commutative diagram
$$
\xymatrix{
Z'\ar[r]^{j'}\ar[d]_h& Y'\ar[d]^g \\
Z\ar[r]^j&Y\ar[r] &  X }
$$
made of closed immersions and such that the square is cartesian. Then for any $(a,\,
m) \in \EE_{Z}(Y) \times \MM_{Y'}(X)$ the relation $h_\flat (g^*(a)\cdot m)=a\cdot
g_\flat(m)$ holds. Analogous formulas hold for classes in $\MM_{Z}(Y) \times
\EE_{Y'}(X)$.

\end{enumerate}

\end{propo}
\qed

\subsection{Arithmetic and relative cohomologies}
\label{Sec Ejemplos Modulos}

We introduce the notion of relative cohomology of a morphism in the context of stable
homotopy theory of schemes and review Holmstrom-Scholbach's construction the
arithmetic counterparts of higher $K$-theory and motivic cohomology. Both
constructions define modules over cohomology theories. We use the theory of monoids
and modules in model categories. This theory can be found written in the context of
motivic homotopy theory in \cite[\S 7]{CD} and a more accessible summary for $\SH$ in
\cite[\S 2.2]{Deglise3}.

\begin{defi}
We say that a ring spectrum, an absolute ring spectrum, a module, an absolute module,
or a morphism is \textbf{strict} if it is defined in the category of symmetric
spectra $\Spt_\Sigma(X)$ and all diagrams commute in $\Spt_\Sigma(X)$ and not just in
$\SH(X)$. Thus, a strict ring spectrum $E$ is a commutative monoid in the category
$\Spt_\Sigma(X)$ and a morphism of strict ring spectra is morphism in
$\Spt_\Sigma(X)$ compatible with the structure morphisms. We denote by $\Mon (X)$ the
category of strict ring spectra with morphism morphisms of strict ring spectra. Let
$E$ be a strict ring spectrum, we denote $E\mod$ the category of strict $E$-modules
with strict morphisms of modules.
\end{defi}

\begin{ejem}
The spectra from Example \ref{Ejem espectros en anillo} representing cohomologies are
strict. Find a reference for the algebraic $K$-theory spectrum $\KGL$ in
\cite[13.3.1]{CD}, for $\HB$ in \cite[14.2.6]{CD}, for mixed Weil theories,
Deligne-Beilinson and absolute Hodge spectra one checks that the product and unit
morphism from \cite[2.1.5]{CD2} for mixed Weil theories and from \cite[1.4.10]{DM}
for Deligne-Beilinson and absolute Hodge spectra satisfy the monoid axioms. Moreover,
the Chern character and the cycle class map are morphism of strict ring spectra
\cite[14.2.16]{CD}. Spitzweck's motivic cohomology spectrum is also strict
(\cite[p.4]{Spitz}).
\end{ejem}

\begin{obser}\label{Obser modulos estrictos}
The categories $\Mon (X)$ and $E\mod$ inherit a model structure from the
$\AA^1$-stable symmetric model structure in $\Spt_\Sigma(X)$. The categories $\Ho
(\Mon (X))$ are well behaved with respect to inverse and direct image as described in
\cite[7.1.11]{CD}. We will use the following fact: let $f\colon Y \to X$ be a
morphism of schemes, $E\in \Mon (Y)$ and $F\in \Mon (X)$, then $f_*E$ in $\SH(X)$ and
$f^*E$ in $\SH(Y)$ are given by strict ring spectra.

Let $\EE$ be a strict absolute ring spectrum. The categories $\Ho (\EE_{X}\mod)$ also
have good functorial properties. Moreover, they are triangulated categories and the
forgetful functor
$$
\Ho (\EE_X\mod)\to \SH(X)
$$
is triangulated.
\end{obser}

\bigskip

Let $\varphi\colon \EE\to \FF$ be a morphism of strict absolute ring spectra and $X$
be a scheme. Denote $\varphi_X\colon \EE_X \to \FF_X$ the morphism of strict ring
spectra in $\Spt_\Sigma(X)$ and $\hofib (\varphi_X)$ the homotopy fiber of
$\varphi_X$. The spectrum $\hofib(\varphi_S)$ in $\Spt_\Sigma(S)$ defines by pullback
an absolute spectrum, which we denote $\hofib (\varphi)$. Recall that the homotopy
fiber fits into a distinguished triangle. In other words, for $X$ an $S$-scheme we
have that
\begin{equation}\label{Ecua triangulo hofib}
\hofib (\varphi_X)\longrightarrow \EE_X
\buildrel{\varphi_X}\over{\longrightarrow}\FF_X\longrightarrow \hofib(\varphi_X)[1].
\end{equation}
Since $\EE_X\mod$ is a model category triangulated and $\varphi$ is a morphism of
absolute spectra the following result is straightforward.

\begin{propo}
Let $\varphi\colon \EE \to \FF$ be a morphism of strict absolute ring spectra. With
above notations, $\hofib (\varphi)$ is a strict absolute $\EE$-module and for every
$S$-scheme $X$ we have $\hofib (\varphi)_X=\hofib (\varphi_X)$.
\end{propo}
\qed

\begin{obser}\label{Obser ideal}
Still with above notations, after a replacement we can assume $\varphi_S$ to be a
fibration and $\FF_S$ to be fibrant so that $\hofib(\varphi_S)$ fits into a cartesian
square
$$
\xymatrix{\hofib (\varphi_S)\ar[d]\ar[r]&\EE_S\ar[d]^{\varphi_S}\\
\mbox{*}\ar[r]&\FF_S.}
$$
Note that the replacement is functorial so we have a commutative diagram
$$
\xymatrix{\hofib(\varphi_{S})\wedge \hofib(\varphi_{S})\ar[r]\ar@{-->}[d] &\EE_{S}
\wedge \EE_{S} \ar[d]^-\mu \ar[r]^{\varphi_S\wedge \varphi_S}&\FF_S\wedge
\FF_S \ar[d]^{\bar \mu}\\
\hofib(\varphi_{S}) \ar[r]&\EE_{S}\ar[r]^\varphi& \FF_S}
$$
Therefore, the groups $\Hom_{\SH(X)} (\11_{X} ,\hofib \varphi_X)$ not only are
modules over $\EE(X)$ but also have an inner product. Note that they do not have a
unit. The distinguished triangle (\ref{Ecua triangulo hofib}) gives raise to a long
exact sequence
$$
\cdots \to \FF^{*-1,*}(X)\to \Hom_{\SH(X)}(\11_X, \hofib (\varphi)(*)[*])\to
\EE^{*,*}(X)\to \FF^{*,*}(X)\to \cdots
$$
where arrows are compatible with products.
\end{obser}

\medskip

We introduce the relative cohomology in the context of motivic homotopy. Let $f\colon
Y\to X$ be a morphism of schemes, then $f_{*}\EE_{Y}$ represents in $\SH(X)$ the
cohomology of $Y$. Indeed,
$$
\EE^{*,*}(Y)=\Hom_{\SH(Y)}(f^*\11_{X}, \EE_{Y}(*)[*])=\Hom_{\SH(X)}(\11_{X},
f_{*}\EE_{X}(*)[*]).
$$
Since $\EE_Y\simeq f^*\EE_X$ we have an adjunction morphism $\EE_X \to f_*f^*\EE_X$.
Let us remark two properties of this morphism.

\begin{propo}\label{Propo f_*E}
Let $\EE$ be an absolute ring spectrum and $f\colon Y\to X$ be a morphism of schemes:
\begin{enumerate}
\item
The spectrum $f_*\EE_Y$ is a ring spectrum. The adjunction $\EE_X\to f_*\EE_Y$ is a
morphism of ring spectra and it induces the inverse image on cohomology $f^*\colon
\EE(X)\to \EE(Y)$.

\item
If in addition $\EE$ is strict, then $f_*\EE_Y$ is also strict and the adjunction map
$\EE_X \to f_*\EE_X$ is represented by a morphism of strict ring spectra.

\end{enumerate}
\end{propo}
\qed

\begin{defi}
Let $\EE$ be a strict absolute ring spectrum, $f\colon Y\to X$ be a morphism of
schemes. Abuse notation and denote $\hofib_\EE(f_X)$ (or simply $\hofib(f_X)$ if it
is clear by the context) the strict $\EE_X$-module defined as the homotopy fiber of
the strict morphism of ring spectra $\EE_X\to f_*\EE_Y$. We define the
\textbf{relative cohomology} of $f$ to be
$$
\EE^{p,q}(f)\coloneqq \Hom_{\SH(X)}(\11_X, \hofib _\EE(f_X)(q)[p]) \quad \mbox{for
}p,q \in \ZZ .
$$
We also denote $\hofib _\EE(f)$ (or simply $\hofib(f)$) the strict absolute
$\EE$-module over $X$ that $\hofib_\EE(f_X)$ defines by pullback.
\end{defi}

\begin{obser}
Consider the above notations. Note that, although $\EE$ is an absolute spectrum, the
spectrum $\hofib (f)$ need not to have good functorial properties. More concretely,
consider a cartesian square
$$
\xymatrix{Y_T\ar[r]\ar[d]_{f_T}&Y\ar[d]^f\\
T\ar[r]^g &X.}
$$
Then the spectrum $g^*f_*\EE_Y$ may not be isomorphic to $f_{T*}\EE_{Y_T}$. In other
words, the family of spectra $f_{T*}\EE_{Y_T}$ for $T\to X$ may not define an
absolute spectrum. Therefore $\hofib(f_T)$ may not be isomorphic to $\hofib
(f)_T=g^*\hofib(f_X)$ so that $\hofib(f)$ may not represent the relative cohomology
of $f_T$.
\end{obser}

\begin{propo}\label{Propo cambio base liso hofib}
Let $\EE$ be a strict absolute spectrum, $f\colon Y\to X$ and $g\colon T\to X$ be two
morphism of schemes. If either $f$ is proper or $g$ is smooth then we have
$$
g^*\hofib (f_X)\simeq \hofib (f_T).
$$
\end{propo}
\demo Denote $g'\colon Y_T\to Y$ and $f_T\colon Y_T\to T$. It is enough to prove that
$g^*f_*\EE_Y\simeq f_{T*}{g'}^*\EE_Y$. The result follows from the smooth and proper
base change in stable homotopy (\cf \cite[1.1.19 and 2.4.50]{CD}).

\qed

\begin{teo}
Let $f\colon Y\to X$ be a morphism of regular schemes and denote $K(f)$ the relative
algebraic $K$-theory (\cf \cite[IV.8.5.3]{Weibel K libro}). Then
$$
\qquad \phantom {a} \qquad K_i(f)=\Hom_{\SH(X)}(\11_X, \hofib_\KGL (f_X)[-i]) \quad
\mbox{for }i \in \ZZ.
$$
\end{teo}
\demo We use notation from \cite{CD}. Recall from \cite[\S 3.2]{CD} that there is
total derived global section functor
$$
\mathrm{R}\Gamma \colon \SH(X)\to \Ho (\Spt_{S^1}).
$$
where $\Spt_{S^1}$ denotes the classic category of $S^1$-spectra of simplicial sets.
Recall from \cite[13.4]{CD} that $K_n(X)=\pi_n(\mathrm{R}\Gamma
(X,\KGL_X))=\Hom_{\SH(X)}(\11_X,\KGL_X)$. Applying the total derived global section
functor to the homotopy fiber sequence
$$
\hofib_\KGL (f_X)\to \KGL_X\to f_*\KGL_Y
$$
we obtain a (classic) homotopy fiber sequence
$$
\mathrm{R}\Gamma (X,\hofib_\KGL (f_X))\to \mathrm{R}\Gamma (X,\KGL_X)\to
\mathrm{R}\Gamma (X,f_*\KGL_Y).
$$
We conclude by recalling that the relative $K$-theory of $f$ is defined in
\cite[IV.8.5.3]{Weibel K libro} as the classic homotopy fiber.

\qed

\begin{ejem}\label{Ejem cohomologia relativa}
Let $\EE$ be a strict absolute ring spectrum. The construction of the relative
cohomology generalizes many concrete situations:
\begin{enumerate}
\item
Let $S=\Spec (k)$ and $p\colon X\to S$ be the structural morphism. Then
$\EE(p)=\widetilde{\EE}(X)$, the reduced cohomology of $X$.

\item
Let $i\colon Z\to X$ be a closed immersion with open complement $j\colon U\to X$.
Recall from \cite[2.3.3]{CD} that we have a distinguished triangle
$$
i_*i^!\EE_X\to\EE_X\to j_*j^*\EE_X\to i_* i^!\EE_X[1].
$$
Therefore $\hofib (j)\simeq i_*i^!\EE_X$ in $\SH(X)$ and $\EE(j)=\EE_Z(X)$. The
product as $\EE$-module is the same of Definition \ref{Defi producto refinado}.

\item
Consider the above notations. Then we also have a distinguished triangle
$$
j_\sharp j^*\EE_X\to \EE_X\to i_*i^*\EE_X\to j_\sharp j^*\EE_X[1]
$$
so that $\hofib (i)\simeq j_\sharp j^*\EE_X$. Although we have not reviewed it, $j_!
j^*\EE_X$ represents, by definition, the cohomology of $U$ with \emph{compact
support} $\EE_\mathrm{c}(U)$ in $\SH(X)$.

\item
Let $i\colon Z\to X$ be a closed immersion and consider the blow-up cartesian square
$$
\xymatrix{P\ar[d]_{\pi'}\ar[r]&B_ZX\ar[d]^\pi\\
Z\ar[r]^i& X.}
$$
It follows from upcoming Corollary \ref{Coro inyectivo} that if in addition $\EE$ is
oriented (\cf \ref{defi orientacion}) then $\EE(\pi)=\EE(P)/ \EE(Z)$. The product as
$\EE(X)$-module is the product through $(i\pi')^*\colon \EE(X)\to \EE(P)$.

\item
Let $R$ be a Dedekind domain and $F$ be its field fractions. Denote $K(R)=K(\Spec
(R))$, $K(F)=K(\Spec (F))$ and $\gamma \colon \Spec(F)\to \Spec (R)$ the localization
morphism. Then $K_2(\gamma)=\coprod_{\p}K_2(R/\p)$ where $\p$ denote prime ideals of
$R$ and $K_1(\gamma)=\coprod_{\p}(R/\p)^\times$ (\cf \cite[III.6.5]{Weibel K libro}).

\end{enumerate}
\end{ejem}

We review the construction of the arithmetic counterparts of $K$-theory and motivic
cohomology of \cite{HS}, which are another example of a homotopy fiber.

Let $A$ be an arithmetic ring (\cf Appendix \ref{Defi cuerpo aritmetico}) and
$S=\Spec (A)$ and $\eta$ its generic point. Recall from Example \ref{Ejem espectros
en anillo} that we have the Deligne-Beilinson cohomology strict ring spectrum
$\EE_{\mathrm{D},\eta}\in \SH (\eta)$, which defines a strict absolute ring spectrum
$\eta_*\EE_{\mathrm{D},\eta}\in \SH(S)$. Recall from Example \ref{Ejem espectros en
anillo} that we have the cycle class map $\cl \colon \mathrm{H}_{\cyrrm{B},\eta} \to
\EE_{\mathrm{D},\eta}$ which induces a strict map
$$
\varphi\colon \mathrm{H}_{\cyrrm{B},S}\to
\eta_*\mathrm{H}_{\cyrrm{B},\eta}\xrightarrow {\eta_*\cl}
\eta_*\EE_{\mathrm{D},\eta}.
$$
Recall that in $\SH(S)_\QQ$ we have $\KGL_{S,\QQ}=\bigoplus_{i\in
\ZZ}\mathrm{H}_{\cyrrm{B},S}(i)[2i]$ (\cf \cite[\S 14]{CD}). We have a strict map
$$
\oplus (\ch_i \circ \varphi_i)\colon \KGL_{S,\QQ}\to
\eta_*\EE_{\mathrm{D},\eta}(i)[2i]
$$
where $\varphi_i\colon \mathrm{H}_{\cyrrm{B},S}(i)[2i]\to
\eta_*\EE_{\mathrm{D},\eta}(i)[2i]$.

\begin{defi}
In above notations, we define the \textbf{arithmetic motivic cohomology} strict
absolute spectrum as $ \widehat {\mathrm{H}}_{\cyrrm{B},S} =\hofib (\varphi ).$ Let
$X$ be a smooth $S$-scheme, we denote the cohomology it defines as
$$
\widehat {H}_\M ^p(X,q)\coloneqq \Hom_{\SH(X)_\QQ}(\11_X,
\widehat{\mathrm{H}}_{\cyrrm{B}, X}(q)[p]) \quad \mbox{for }p,q \in \ZZ .
$$

Analogously, we define the \textbf{arithmetic $K$-theory} strict absolute spectrum as
$ \hatKGL_{S,\QQ} =\hofib (\oplus (\ch_i \circ \varphi_i))$. Note that the
periodicity of the $K$-theory makes $\hatKGL$ also periodic. Let $X$ be a smooth
$S$-scheme, we denote the cohomology it defines as
$$
\hatKH _i (X)_\QQ \coloneqq \Hom_{\SH(X)_\QQ}(\11_X, \hatKGL_\QQ [-i]) \quad
\mbox{for } i \in \ZZ.
$$
\end{defi}

\begin{obser}
Note that our definition is written differently from \cite{HS}, where they considered
the spectra
$$
\hofib(\mathrm{H}_{\cyrrm{B},S}\buildrel{\mathrm{id}\wedge
1_\mathrm{D}}\over{\longrightarrow} \mathrm{H}_{\cyrrm{B},S} \wedge \eta_*
\EE_{\mathrm{D}, \eta})\mbox{ and } \hofib (\KGL_{S} \buildrel{\mathrm{id}\wedge
1_\mathrm{D}}\over{\longrightarrow} \KGL_S \wedge \eta_* \EE_{\mathrm{D}, \eta}).
$$
Recall that in $\SH(S)_\QQ$ we have $\mathrm{H}_{\cyrrm{B},S} \wedge \eta_*
\EE_{\mathrm{D}, \eta}\simeq \eta_*\EE_{\mathrm{D}, \eta}$ (\cf \cite[14.2.8]{CD}) so
that both definitions agree.
\end{obser}

\begin{obser}\label{Caracter de Chern aritmetico}
By construction, both $\hatHm$ and $\hatKGL_\QQ$ are strict absolute $\HB$ and
$\KGL_\QQ$ modules respectively, but not rings. Nevertheless, the Chern character
$\ch\colon \KGL_\QQ \to \HB$ induces an \emph{arithmetic Chern character} $
\hatch\colon \hatKGL _\QQ\to \hatHm $ and the square
$$
\xymatrix{\hatKGL _\QQ \ar[r] \ar[d]_{\hatch}&\KGL_\QQ \ar[d]^{\ch}\\
\hatHB \ar[r]&\ \HB}
$$
commutes.
\end{obser}

\subsection{Orientations}

We review the theory of orientations (i.e., Chern classes) for spectra. As in the
classical case, they are determined by the first Chern class of the tautological line
bundle of projective spaces.

Recall that $\11_S(1) = \mathrm{coker} (\Sigma^\infty S\to
\Sigma^{\infty}\PP^1)[-2]$. For any ring spectrum $\EE$ with unit $\eta\colon S\to
\EE_S $ there is a canonical class in $\EE^{2,1}(\PP^1)$ defined as the morphism
$$
\PP ^1 \to \11_S(1)[2]=S\wedge \11_S(1)[2] \xrightarrow{\eta\wedge \mathrm{Id}} \EE_S
(1)[2].
$$
By abuse of notation we still denote it $\eta\in \EE ^{2,1}(\PP^1)$.

The definition of $\EE$-cohomology may be extended to general spectra. In particular,
recall that the infinite projective space is defined to be $\PP^\infty_X=\varinjlim
\PP^n_X$ and we denote
$$
\EE^{p,q} (\PP^\infty_X)= \Hom _{\SH(X)}(\PP^\infty_X,\EE_X (q)[p]).
$$

\begin{defi}\label{defi orientacion}
We define an \textbf{orientation} on an absolute ring spectrum $\EE$ to be a class
$c_{1} \in \EE ^{2,1} (\PP ^\infty)$ such that for $i_1\colon \PP^1 \hookrightarrow
\PP^\infty$ satisfies $i_1 ^* (c_{1})= \eta$. We also say that $\EE$ is
\textbf{oriented}.

Let $X$ be a scheme and $\V$ be a locally free $\O_X$-module. We call the
\textbf{vector bundle} given by $\V$ to the scheme $V=\mathbf{Spec}\, _X (S^\bullet
\V^*) \to X$ and the \textbf{projective bundle} given by $\V$ to the scheme $\PP
(V)=\mathbf{Proj}\, _X (S^\bullet \V^*)\to X$.

Let $B\mathbb{G}_m$ be the classifying space for $\mathbb{G}_m$-torsors
(\cite[4.1.16]{MV}), due to \cite[4.3.7]{MV}  we have a natural map
$$
\Pic (X)  \to \Hom_{\HH_\bullet (X)}(X_+,B\GG_m)\simeq \Hom _{\HH_\bullet (X)}
(X_+,\PP^\infty) \to \Hom_{\SH(X)}(\11_X,\PP^\infty _X)
$$
so that any line bundle $L\in \Pic (X)$ defines a morphism $f\colon X\to \PP ^\infty
_X$ in $\SH (X)$.

Let $(\EE, c_1)$ be an oriented absolute ring spectrum and denote
$p\colon\PP^\infty_X\to \PP^\infty$. For any line bundle $L$ we have
$$
\begin{matrix}
\EE^{2,1}(\PP^\infty_X) &\xrightarrow {f^*} &\EE ^{2,1}(X)\\
p^*c_{1} & \mapsto & c_1( L)
\end{matrix}
$$
and we say that $c_1(L)\coloneqq f^*p^*c_1$ is \textbf{the first Chern class} of $L$.
\end{defi}

\begin{ejem}\label{Ejem orientados}
Every example of Example \ref{Ejem espectros en anillo} representing a cohomology is
oriented. We review the references: Mixed Weil theories are oriented in
\cite[2.2.8]{CD2}, the algebraic $K$-theory $\KGL$ and Beilinson's motivic cohomology
$\HB$ in \cite[13.2.2]{CD} and \cite[14.1.5]{CD} respectively, algebraic cobordism
$\MGL$ in \cite[1.4]{Panin3}. In particular, the orientation of $K$-theory is given
by $c_1^{\KGL}(L)=1-[L^*]$. Spitzweck's motivic cohomology spectrum
$\mathrm{H}_\Lambda$ is oriented in \cite[11.1]{Spitz}. In the Appendix \ref{Hodge
absoluto} we give an orientation for the absolute Hodge spectrum and the
Deligne-Beilinson is done in \cite[3.6]{HS}. Finally, every cohomology considered in
\cite{DM} is represented by an oriented spectrum (\emph{cf.} \cite[1.4.11.(1) and
2.1.2.(1)]{DM}). In particular, Besser's absolute rigid syntomic spectrum is
oriented.
\end{ejem}

\begin{obser}
Let $\varphi \colon \EE\to \FF $ be a morphism of absolute ring spectra and let
$c_1\in \EE^{2,1}(\PP^\infty)$ be an orientation on $\EE$. Since $\varphi$ is a
morphism of rings, it maps the unit $\11_S\to \EE_S$ onto the unit $\11_S \to \FF_S$.
We conclude that the element $\varphi_{\PP^\infty}(c_1)\in \FF^{2,1}(\PP^\infty)$ is
an orientation on $\FF$.
\end{obser}

\begin{obser}
Denote $i_n\colon \PP^n\to \PP^\infty$ since $i_n^*(\O(-1))\simeq \O_{\PP^n(-1)}$ we
have that $i_n^*(c_1)=c_1(\O_{\PP^n}(-1))$ and we write $c_1 = c_1
(\O_{\PP^\infty}(-1))$.
\end{obser}

To fix notations, we recall  the theory of Chern classes in the context of spectra.
Proof of the following result in the context of stable homotopy theory may be found
in \cite[2.1.13 and 2.1.22]{Deglise}.

\begin{teo} [Projective bundle]\label{Teo Hirsch Leray}
Let $V\to X$ be a vector bundle of rank $(n+1)$, $\EE$ an oriented absolute ring
spectrum and $x=c_1(\O_{\PP (V)} (-1))$. There is a canonical isomorphism
$$
\begin{matrix}
\bigoplus_{i=0}^n \EE ^{*-2i,*-i} (X) &\buildrel{\sim}\over{\longrightarrow } &
\EE ^{*,*} (\PP (V))\\
\ \ \ \ \ \ (a_0, \ldots , a_n) & \mapsto & \sum _i \pi ^*  (a_i) x^i.
\end{matrix}
$$
\end{teo}

\qed

\begin{defi}
Let $V\to X$ be a vector bundle of rank $n$. We define the \textbf{$i$-th Chern
classes} $c_i(V)\in \EE^{2i,i}(X)$ for $i=1 , \ \ldots , \ n$ as the unique ones
satisfying
$$
c_1(\O_{\PP(V)}(-1))^n+\sum _{i=1}^{n} (-1)^i c_i (V) c_1( \O_{\PP (V)}(-1))
^{n-i}=0.
$$
\end{defi}

\begin{parra}\label{serie f}
\emph{Formal (abelian) group laws} $F(x,y)$ are certain type of series
(\cf \cite{Adams}). In
particular, they satisfy the property that any formal group law $F(x,y)$ is of the
form
$$
F(x,y)=x+y+xyf(x,y)
$$
for $f(x,y)$ a symmetric series. A formal group law is called \textbf{additive} if
$F(x,y)=x+y$. The following is a classic result. See for example
\cite[3.7]{Deglise2}.
\end{parra}

\begin{teo}\label{Ley de Grupo Formal de c1}
Let $\EE$ be an oriented absolute ring spectrum. There exists a formal group law
$F(x,y)\in \EE^{**}(S)[[x,y]]$ such that
$$
c_1(L_1 \otimes L_2)=F(c_1(L_1), c_1(L_2))
$$
for any line bundles $L_1$, $L_2$ over $X$. In addition, Chern classes are nilpotent.
\end{teo}
\qed

See \cite[6.2]{NSO} for a proof of the following result in our context.

\begin{propo}
Let $(\EE,c_1)$ be an oriented absolute ring spectrum and denote
$\EE(\PP^\infty)=\prod \EE^{p,q}(\PP^\infty)$, then
$$
\EE(\PP^\infty)=\EE(S)[[c_1]].
$$
\end{propo}
\qed

\begin{ejem}
Every cohomology from Example \ref{Ejem espectros en anillo} apart from $K$-theory
and algebraic cobordism have additive formal group laws. That is to say: motivic
cohomology, cohomologies coming from mixed Weil theories, real absolute Hodge and
Deligne-Beilinson cohomology and Besser's rigid syntomic cohomology have additive
formal group laws.
\end{ejem}

Standard arguments yield the following classic formula:

\begin{teo}[Whitney sum] \label{Cartan-Weil}
Let $\EE$ be an oriented absolute ring spectrum and $ 0 \to V' \to V \to V'' \to 0 $
be a short exact sequence of vector bundles, then we have
$$
c_k(V)=\sum _{i+j=k}c_i(V')c_j(V'') \ \ i,j,k \in \NN.
$$
\end{teo}

\qed

\begin{propo}\label{serie c1}
Let $(\EE, c_1)$ be an oriented absolute spectrum and $c_1^{\mathrm{new}}$ be another
orientation. Then there exists $G(t)\in \EE (S)[[t]]$ with leading coefficient 1 such
that for any line bundle $L$ we have
$$
c_1^{\mathrm{new}}(L)=G(c_1(L))c_1(L).
$$
\end{propo}
\demo Since $\PP^{\infty}\simeq \mathrm{B}\GG_m$, the classifying space for line
bundles, it is enough to check the formula for $x=c_1(\O_{\PP^\infty}(-1))$. Recall
that $\EE(\PP^{\infty})=\EE (S)[[x]]$ and therefore we have
$c_1^{\mathrm{new}}(\O_{\PP^\infty}(-1))= a_0+a_1x+\ldots =G(x)x $.
Finally, both classes satisfy
$i_1^*(c_1(\O_{\PP^{\infty}}))=i_1^*(c_1^{\mathrm{new}}(\O_{\PP^{\infty}}))=\eta\in
\EE ^{2,1}(\PP^1)$ and we conclude that $a_0=0$ and $a_1=1$.

\qed

\medskip

\begin{parra}\label{Tod}
We recall in our context the multiplicative extensions of a series. Let $\EE$ be an
absolute oriented ring spectrum and $F(t)\in \EE(S)[[t]]^\times$ be an invertible
series. Let $L\to X$ be a line bundle, we set
$$
F(L)\coloneqq F(c_1(L))\in \EE(X)
$$
(it is well defined since Chern classes are nilpotent). Let $V\to X$ be a vector
bundle, by the splitting principle we can assume $V=L_1+\ldots + L_n$. We put
$$
F_\times (V)\coloneqq F(L_1)\cdot \ldots \cdot F(L_n) \in \EE(X)
$$
where we consider $F(L_1)\cdot \ldots \cdot F(L_2)$ as a power series in the
symmetric elemental functions of $c_1(L_1),\ldots , c_1(L_n)$, which are just the
Chern classes of $E$. We call $F_\times (T)$ the \emph{multiplicative extension} of
$F$.

\end{parra}

\begin{coro}\label{Tod cn}
Let $(\EE, c_1)$ be an oriented absolute spectrum and $c_1^\mathrm{new}$ be another
orientation on $\EE$ and denote $G(t)\in \EE(S)[[t]]^\times$ the series such that
$c_1^\new(L)=G(c_1(L))c_1(L)$. Then
$$
c_n^\new (V)=G_\times (V)c_n(V)
$$
for $V$ a rank $n$ vector bundle.
\end{coro}

\qed

\subsection{Chern class with support}

\begin{defi}
Let $X$ be a scheme and $U=(X-Z)\to X$ be an open subscheme. We call a \textbf{pseudo
divisor} (trivialized on $U$) a pair $(\L,u)$ consisting of an invertible sheaf $\L$
over $X$ (in the Zariski topology) and a trivialization $u\colon\O| _U\xrightarrow
\sim \L | _U$. We denote $\Pic_Z (X)$ the group of isomorphism classes of pseudo
divisors.
\end{defi}

The same arguments of \cite[\S 4.1.3]{MV} give the following remark:

\begin{obser}
Let $Z\to X$ be a closed immersion of codimension 1 and $U$ be its open complement.
Then there is a map
$$
\Pic_Z (X) \to \Hom_{\HH_\bullet^s (X)}(X/U , \PP^\infty _{X}).
$$
\end{obser}

It follows from above identification that there is a map \goodbreak
$$
\mathrm{Pic}_Z(X)\to \Hom_{\HH_\bullet(X)}(X/U,\PP_{X} ^\infty)\to \Hom_{\SH
(X)}(\Sigma^\infty X/U, \PP_X ^\infty) \to
$$
$$
\buildrel{(c_{1})_*}\over{\to}\Hom_{\SH  (X)}(\Sigma^\infty X/U, \EE_X (1) [2])=\EE
^{2,1}_Z(X).
$$
Finally, for an oriented absolute spectrum $ \EE$ and any line bundle $L$ we have
$$
\begin{matrix}
\Pic_Z(X) &\xrightarrow {\varphi_{L,u} ^*} &\EE_Z  ^{2,1}(X)\\
(L,u) & \mapsto & c_1^Z( L,u)
\end{matrix}
$$
and we say that $c_1^Z(L,u)\coloneqq \varphi_{L,u}^*(L,u)$ is the first Chern class
of $L$ with support on $Z$. We omit the reference to the trivialization when no
confusion is possible. The next statement follows from the definition:

\begin{propo} \label{Propo funtorialidad clase Chern refi}
Let $f\colon X'\to X$ be a morphism of schemes and $(\L,u)$ be a pseudo divisor over
the open $X-Z$, then
$$
f^* c_1 ^Z(L)=c_1 ^{f^{-1}Z}(f^*L).
$$
\end{propo}

\qed

\begin{propo}\label{Propo serie f soporte}
Let $\L_1$ and $\L_2$ be two invertible sheaves over $X$ and $u_i\colon  \O|_U\to
\L_i|_U$, $i=1,2$, trivializations. Denote  $\L =\L _1 \otimes \L _2$ and $u=u_1
\otimes u_2$. Let $F(x,y)\in \EE (S)[[x,y]]$ be the formal group law of $c_1$ given
by Theorem \ref{Ley de Grupo Formal de c1}. Then
$$
c_1 ^Z(L)=F(c_1 ^Z (L_1),c_1 ^Z(L_2)).
$$
\end{propo}
\demo The pseudo divisors $(L_1, u_1)$, $(L_2,u_2)$ and $(L,u)$ correspond to
morphism $f_1$, $f_2$, $f\colon X/U \to \PP_X^\infty$ respectively. Denote the Segre
embedding $\sigma\colon\PP_X^\infty \times \PP _X^\infty \to \PP_X^\infty$. By
construction, the diagram
$$
\xymatrix{ X/U \ar[r]^-{f_1\times f_2} \ar[dr]_f & \PP_X ^\infty \times
\PP _X ^\infty \ar[d]^\sigma\\
& \PP _X ^\infty }
$$
commutes and, taking $\EE$-cohomology, we have the commutative diagram
$$
\vcenter{\xymatrix{ \EE(\PP_X ^\infty) \ar[r]^-{\sigma^*} \ar[dr]_-{f^*} &
\EE(\PP_X ^\infty \times\PP_X ^\infty)\ar[d]^-{(f_1\times f_2)^*}\\
& \EE_Y(X).}}
$$

Recall that $\EE(\PP^\infty)=\EE(S)[[t]]$ and $\EE(\PP^\infty \times
\PP^{\infty})=\EE(S)[[u,v]]$ where $t=c_1(\O_{\PP^\infty}(-1))$,
$u=c_1(p_1^*\O_{\PP^\infty}(-1))$, $v=c_1(p_2^*\O_{\PP^\infty}(-1))$ and $p_1$, $p_2$
are the canonical projections. With this notations, the Segre morphism maps $t\mapsto
F(u,v)$ where $F$ is the formal group law of the orientation of $\EE$. We conclude by
the commutativity last diagram.

\qed

\section{Gysin morphism}

\subsection{Regular immersions}\label{Sec inmersion regular}

We construct the Gysin morphism for a regular immersion following Gabber's ideas for
\'{e}tale cohomology (see \cite{Fujiwara} and \cite{Riou}). More concretely, we lift the
approach of Riou in \cite{Riou} for \'{e}tale cohomology to the motivic homotopy setting.
The main difference lies in the fact that \'{e}tale cohomology has additive Chern
classes, that is $c_1(L_1\otimes L_2)=c_1(L_1)+c_1(L_2)$ for $L_1$ and $L_2$ line
bundles (\cf Paragraph \ref{serie f}), while in the general setting Chern classes are
not additive (\cf Example \ref{Ejem orientados}).

Gabber's method reduces the case of general codimension to that of codimension one by
means of describing the cohomology of the blow-up (\emph{cf.} Corollary \ref{Coro
inyectivo}). In order to prove the needed functoriality properties the versatile
context of modified  blow-up is used.

\begin{defi}
Let $i\colon Z\to X$ be a regular immersion of codimension 1. The sheaf
$\mathcal{I}_Z=\L_{-Z} $ is locally principal and it has a natural trivialization on
$X-Y$, and so does its dual. We define the \textbf{refined fundamental class} (of $Z$
in $X$) to be
$$
\bar \eta_Z ^X \coloneqq c_1 ^Z (\I_Z^*)=c_1^Z(L_Z)\in \EE_Z^{2,1} (X)
$$
and the \textbf{fundamental class} to be $\eta _Z ^X \coloneqq  c_1 (L_Z)=i_\flat
(c_1^Z(L_Z)) \in \EE ^{2,1}(X)$. We define the \textbf{refined Gysin morphism} as
$$
\begin{matrix}
\p_i \colon & \EE ^{*,*}(Z)& \longrightarrow & \EE^{*+2,*+1}_{Z}(X)
\\ & a & \mapsto & a\cdot \bar \eta _Z ^X
\end{matrix}
$$
and the \textbf{Gysin morphism} as $i_*\colon\EE ^{*,*}(Z) \longrightarrow
\EE^{*+2,*+1}(X)$, $a  \mapsto  i_\flat(a\cdot \bar \eta _Z ^X)$.

More generally, let $(\L , u\colon\O |_U\xrightarrow\sim \L | _U)$ be a pseudo
divisor where $U=X-Z$. We define the \textbf{refined Gysin morphism (given by $(\L,
u)$)} as
$$
\begin{matrix}
\p_{\L}\colon& \EE ^{*,*}(Z)& \longrightarrow & \EE^{*+2,*+1}_{Z}(X)
\\ & a & \mapsto & a\cdot c_1 ^Z (L).
\end{matrix}
$$
\end{defi}

\begin{obser}\label{Coro c1 soporte}
Let $Z\xrightarrow i X$ be a regular immersion of codimension 1. Due to Proposition
\ref{Propo propiedades producto} it is easy to check that
$$
c_1^Z (L_Z)\cdot c_1 (N_{Z/X})=c_1^Z (L_Z )^2
$$
where $c_1 ^Z(L_Z)\in \EE_Z ^{2,1}(X)$, $N_{Z/X}=\Spec (S^\bullet (\I_Z/\I_Z^2) ^*)$
and $c_1(N_{Z/X})\in \EE^{2,1} (Z)$.
\end{obser}

We now define a refined fundamental class for any closed subscheme $Z\to X$ and any
epimorphism $\F^* \to \I_Z/\I_Z^2$ where $\F$ is a locally free $\O_Z$-module. This
more general context, due to Gabber, is the suitable one to prove the basic
properties of the Gysin morphism.

\begin{defi}
Let $Z\to X$ be a closed immersion of strictly positive codimension defined by a
sheaf of ideals $\I_Z$. Let $\F ^*\to \I_Z /\I_Z ^2$ be an epimorphism of $\O
_Z$-modules where $\F$ is locally free and consider the $\O_X$-graded algebra
$\mathcal{A}=\bigoplus \mathcal{A}_n$ defined on each degree as the fibre product
$$
\begin{matrix}
\mathcal{A}_n & \longrightarrow &\I_Z ^n \\
\downarrow & & \downarrow  \\
S^n \F^*& \longrightarrow & \I_Z^n/ \I_Z ^{n+1}.
\end{matrix}
$$
We define the \textbf{modified blow-up} as the projective scheme $B_{Z,\F} X\coloneqq
\mathbf{Proj}_X (\mathcal{A})$.
\end{defi}

See \cite[2.2.1.3; 2.2.1.4 and 2.1.5]{Riou} for a proof of the following properties
of the modified blow-up:

\begin{propo}
Let $\pi\colon B_{Z,\F} X \to X$ be a modified blowing-up:
\begin{enumerate}

\item
The epimorphism $\mathcal{A}\to \bigoplus \I^n_Z/\I_Z^{n+1}$ defines a closed
embedding $B_{Z} X\to B_{Z,\F} X$.

\item
If $Z\to X$ is a regular immersion and $\F^*=\I_Z/\I_Z^2$ then $B_{Z,\F} X=B_{Z} X$
is the classic blow-up $B_{Z} X$.

\item
Denote $U=X-Z$, then $B_{Z,\F} X|_{\pi ^{-1} (U)}\simeq U$.

\item
$\pi^{-1}(Z)=\mathbb{P}(F)= \mathbf{Proj}_Z ( S^\bullet \F)$.

\item
For any morphism $p\colon X'\longrightarrow X$ there is a canonical morphism
$$
B_{Z',p^* \F} X'\longrightarrow B_{Z,\F} X \times _X X'
$$
which is a nil-immersion.

\end{enumerate}
\end{propo}

\qed

\begin{propo}\label{Explosion}
Let $i\colon Z \to X$ be a closed immersion and let $\F^* \to \I_Z /\I_Z ^2$ be an
epimorphism of $\O _Z$-modules where $\F$ is locally free. Let $B=B_{Z,\F}X$ be the
modified blowing up, $\pi\colon B\to X$ the canonical morphism and $P=\PP(F)$ the
exceptional divisor. Then for any $q\in \ZZ$ there is a long exact sequence
$$
\cdots \to \EE _Z ^{p,q}(X)\to \EE ^{p,q}_P (B) \oplus \EE^{p,q} (Z) \to \EE ^{p,q}
(P)\to  \EE _Z ^{p+1,q}(X)\to\cdots .
$$
The same long exact sequence holds without support.
\end{propo}

\demo First recall that we call a \emph{cdh}-distinguished square any cartesian
square
$$
\xymatrix{
Z' \ar[r]^{i'} \ar[d]_{\pi'} &X' \ar[d]^\pi \\
Z \ar[r]^{i}& X }
$$
such that $i$ is a closed immersion, $\pi$ is proper and defines an isomorphism
$\pi^{-1}(X-Z)\simeq X-Z$. Following \cite[3.3.8]{CD} any \emph{cdh}-distinguished
square gives homotopy bicartesian squares in $\SH(X)$. In particular, for any
absolute spectrum $\EE$ the diagram in $\SH(X)$
$$
\xymatrix{
\EE_X \ar[r]\ar[d] &\pi_* \EE _{X'}\ar[d] \\
i_* \EE_Z \ar[r] & i' _* \pi_* \EE_{Z'} }
$$
is homotopy bicartesian so that there is a distinguished triangle
$$
\EE_X \to \pi_* \EE _{X'} \oplus i_* \EE_Z \to  i _* \pi '_* \EE_{Z'} \to \EE_X[1].
$$
Applying the functor $\Hom _{\SH(X)}(i_*\11_Z(-q), \raya )$ (or $\Hom
_{\SH(X)}(\11_X(-q), \raya )$ for the case without support) to this triangle in the
case $X'=B=B_{Z,\F} X$ and $Z'=P=\mathbb{P}(F)$ we conclude.

\qed

\begin{coro}\label{Coro inyectivo}
With the preceding notations let $p,\, q\in \ZZ$, we have a split short exact
sequence
$$
0\longrightarrow \EE _Z ^{p,q}(X)\longrightarrow \EE ^{p,q}_P (B)
\buildrel{s}\over{\leftrightarrows} \EE ^{2n,q} (P)/\EE^{p,q} (Z) \longrightarrow 0.
$$
The same exact sequence holds without support.
\end{coro}
\demo The preceding long exact sequence may be rewritten as
$$
\cdots \to \EE _Z ^{p,q}(X)\xrightarrow{\pi^*} \EE ^{p,q}_P(B) \xrightarrow{i^*} \EE
^{p,q} (P)/\EE^{p,q} (Z)\to \cdots.
$$
Denote $x= c_1(\O_{P}(-1))$, by the projective bundle theorem we have
$$
\EE ^{p,q} (P)/\EE^{p,q} (Z)= \EE^{p-2,q-1} (Z)x \oplus \cdots \oplus \EE^{p-2(n
-1),q-n+1} (Z) x^{n-1}.
$$
Since $i^*(c_1 ^P (\O_{B}(-1)))=c_1(\O_P(-1))=x$ we define the section $s$ by giving
its value at the generators $ s(x^i)=(c_1 ^P (\O_{B}(-1)))^i $ and linearity.

\qed

\bigskip

\begin{parra}
Consider once again the notations of Proposition \ref{Explosion}. Recall that $Z$ a
is closed subscheme of $X$, $B\to X$ is a modified blowing-up of $X$ over $Z$ and
$P=\pi^{-1}(Z)$. We now construct a distinguished class in $\EE_Z(X)$ to define the
Gysin morphism.

Although $P$ is not in general of codimension 1, the invertible sheaf $\L
_P=\O_B(-1)$ has a canonical trivialization on $B-P$. Therefore we consider the
refined Gysin morphism $\p_{ \O_B(-1)}$ and the diagram
$$
\xymatrix{ 0 \ar[r] &\EE _Z ^{2n,n}(X)\ar[r]^{\pi^*} &\EE ^{2n,n}_P(B)\ar[r]^-{i^*} &
\EE ^{2n,n} (P)/\EE^{2n,*} (Z) \ar[r] & 0 \\
&  & \EE ^{2n-2,n-1} (P)\ar[u]^{\p _{\O_B(-1)}}& & }
$$
where $n=\rk  \F$. Since $\Sigma_0 ^{n-1} (-1)^{n+1+i} c_i(F) x^{n-i}=c_n (F)=0$ in
$\EE (P)/\EE (Z)$, we define
$$
\mathit{Cl} _{Z,\F} ^X\coloneqq \Sigma_0 ^{n-1} (-1)^{n+1+i} c_i (F) x^{n-i-1}\in \EE
^{2n-2,n-1}(P).
$$
Note that
\begin{align*}
i^*\p_{\O_B(-1)}(\mathit{Cl} _{Z,\F} ^X)&=i^*(c_1^P
(\O_B(-1))(\Sigma_0 ^{n-1} (-1)^{n+1+i} c_i(F) x ^{n-i-1}))\\
&=\Sigma_0 ^{n-1} (-1)^{n+1+i} c_i(F) x^{n-i}=0.
\end{align*}
Therefore, there exists a unique class $\bar \eta _{Z,\F}^X\in \EE^{2n,n}_Z(X)$ such
that $\pi^*{\bar \eta_{Z,\F}^X}=\p_{\O_B(-1)}(\mathit{Cl} _{Z,\F} ^X)$.
\end{parra}

\begin{defi}\label{Definicion Gysin}
Let $Z\to  X$ be a closed subscheme and let $\F^*\to  \I_Z/\I_Z^2$ be an epimorphism
of $\O_Z$-modules where $\F$ is locally free of rank $n$. With the preceding
notations, we define the \textbf{refined fundamental class of $Z$ in $X$ modified by
$\F$} to be the unique class $\bar \eta_{Z, \F} ^X$ in $\EE^ {2n,n} _Z (X)$ such that
$\p_{\O_B(-1)} (\mathit{Cl}_{Z,\F} ^X)=\pi ^* (\bar \eta _{Z,\F} ^X)$.

In the case $i\colon Z\to  X$ is a regular immersion of codimension $n$ and
$\F^*=\I_Z/\I_Z^2$ we call this class the \textbf{refined fundamental class of $Z$ in
$X$} and we denote it $\bar \eta_Z^X\in \EE _Z ^{2n,n}(X)$. We define the
\textbf{refined Gysin morphism} to be
$$
\begin{matrix}
\p_i\colon &\EE ^{*,*}(Z)&\to  &\EE ^{*+2n,*+n}_Z(X)\\
&a & \mapsto & a\cdot \bar \eta_Z ^X
\end{matrix}
$$
and the \textbf{Gysin morphism} to be $i_*\colon \EE ^{*,*}(Z) \longrightarrow
\EE^{*+2n,*+n}(X)$, $a  \mapsto  i_\flat(a\cdot \bar \eta _Z ^X)$.
\end{defi}

\begin{propo}\label{Propo propiedades Gysin}
Let $\EE$ be an absolute ring spectrum and $i\colon Z \to  X$ be a regular immersion
of codimension $n$:

\begin{enumerate}

\item
Denote $N_{Z/X}=\Spec (S^\bullet \I _Z /\I_Z ^2)$ the normal bundle, then
$$
i^*  \eta _Z ^X =c_n(N_{Z/X}) \in \EE^{2n,n} (Z).
$$
In particular, let $V\to  X$ be a rank $n$ vector bundle and $s_0\colon X\rightarrow
V$ be the zero section, then $ s_0 ^*\eta_X ^{V}=c_n(V)\in \EE^{2n,n}(X) $.

\item
Projection formula: The Gysin morphism $i_*$ is a morphism of $\EE (X)$-modules. In
other words,
$$
 a \cdot i_* ( b)= i_* (i^* ( a) \cdot  b) \ \forall \  a \in \EE
(X) \ , \  b \in \EE (Z).
$$

\item
Let $r\colon X\to   Z $ be a retraction of $i$, then $i_*$ is a morphism of $\EE
(Z)$-modules with respect to $r^*$. In other words,
$$
i_*( a\cdot  b )=r^*( a)\cdot i_*( b) \ \ \forall \  a \in \EE(Z).
$$
In particular, $i_*( a)=r^*( a)\cdot  \eta_Z^X$.

\end{enumerate}
\end{propo}
\demo The first point is a consequence of the definition. The projection formula
follows since
$$
i_*(i^*( a)\cdot b)=i_\flat (i^*( a)\cdot b\cdot\bar \eta_Z^X)=  a \cdot i_\flat( b
\cdot\bar \eta_Z^X)= a \cdot i_*( b),
$$
where we have used point 6 of Proposition \ref{Propo propiedades producto}.

For the last point follows note
$$
r^*( a)\cdot i_*( b)= r^*( a)\cdot i_\flat( b\cdot \bar \eta_Z^X)= i_\flat( r^*(
a)\cdot  b\cdot \bar \eta_Z^X)$$
$$
=i_\flat(  a\cdot  b\cdot \bar \eta_Z^X) =i_*( a\cdot  b ),
$$
where we have used point 4 and 6 of Proposition \ref{Propo propiedades producto}.

\qed

\begin{propo} \label{Funtorialidad clase fund}
The refined fundamental class is stable under base change. In other words, let
$p\colon X'\longrightarrow X$ be a morphism of schemes, $Z\to  X$ a closed subscheme
and $\F^*\to  \I_Z/\I_Z ^2$ an epimorphism of $\O _Z$-modules where $\F$ is locally
free, then
$$
p ^* (\bar \eta _{Z,\F}^X)=\bar \eta_{p^{-1}(Z),p^*(\F)} ^{X'}.
$$
\end{propo}
\demo It is enough to check
$$
p^*(\mathit{Cl} _{Z,\F}^X)=\mathit{Cl}_{p^{-1}(Z),p^* \F}^{X '} \in \EE
^{2r-2,r-1}(\mathbb{P}(p^* \F))
$$
where $r=\rk \F$. This follows from the fact that Chern classes are functorial and
the induced morphism $\bar p\colon  B_{p^{-1}Z,p^* \F}X' \to   B_{Z,\F}X$ satisfies
$$
p^*\O_{B_{p^{-1}Z,p^* \F}X'} (-1)=\O_{B_{Z,\F}X} (-1).
$$

\qed

\begin{lema}\label{Lema epimorfismo}
With the preceding notations, let $\F'^*\longrightarrow \F^*$ be an epimorphism of
locally free $\O_Z$-modules of constant rank $r'$, $r$ respectively. Denote $\K^*$ be
the kernel and $F$, $F'$, and $K$ the vector bundles they define. We have the
relation
$$
\bar \eta _{Z,\F'}^X= c_{r'-r}(K)\bar \eta _{Z,\F}^X.
$$
\end{lema}
\demo We check that due to the splitting principle we are reduced to prove the case
where $r'=r+1$.\footnote{Let us remark that this result can be proven without relying
in the splitting principle as shown by Riou in \cite{Riou}.} Indeed, we can assume
the epimorphisms comes is the composition from a sequence of epimorphisms
$\F_1^*\to\F_2^*\to \F_3^*\cdots...$ each one with kernel of rank 1. If the
propositions holds for $r'=r+1$ then
$$
\bar \eta _{Z,\F_1}^X= c_{1}(K_1)\bar \eta _{Z,\F_2}^X=c_{1}(K_1)c_{1}(K_2)\bar \eta
_{Z,\F_3}^X.
$$
Since $c_{1}(K_1)c_{1}(K_2)=c_2(K')$ where $K'$ is the kernel of $F_1\to F_3$. We
conclude repeating this process.

For the case $r=r+1$ it is enough to check that
$$
j^* \mathit{Cl} _{Z,\F'}^X= c_{1}(K)\mathit{Cl} _{Z,\F}^X.
$$

By construction there is short exact sequence $0\to  F \to F' \to K \to 0$ of vector
bundles. Recall from Theorem \ref{Cartan-Weil} that $c_i(F')=c_i(F)+c_1(K)c_{i-1}(F)$
and therefore
\begin{align*}
j^*\mathit{Cl} _{Z,\F'} ^X&=j^*((-1)^{r'+1}\Sigma_0 ^{r'-1} (-1)^{i} c_i (F')
x ^{r'-1-i})\\
& = (-1)^{r}\Sigma_0 ^{r}(-1)^{i} (c_i(F)+c_1(K)c_{i-1}(F))x^{r-i}\\
&=(-1)^{r}\Sigma_0 ^{r}(-1)^{i} c_i(F)x^{r-i}\ + (-1)^{r+1} c_1(K)
\Sigma _0 ^{r-1}(-1)^{j} c_j(F)x^{r-1-j}\\
&=c_1(K)\mathit{Cl} _{Z,\F}^X.
\end{align*}

\qed

\medskip

In Gabber's versatile context of modified blow-up the so called \emph{key formula}
(\emph{cf.} \cite[6.7]{Fulton}) and the more general \emph{excess intersection
formula} are a direct consequence of the definition by Proposition \ref{Funtorialidad
clase fund} and Lemma \ref{Lema epimorfismo}.

\begin{coro}[Excess intersection formula]\label{Coro formula exceso}
Consider the cartesian square
$$
\xymatrix{ P \ar[r]^j \ar[d]_{\pi '} &X'\ar[d]^\pi \\
Z\ar[r]^i & X}
$$
where both $i$ and $j$ are regular immersions of codimension $n$ and $m$
respectively. If $K=\pi'^* N_{Z/X}/N_{P/X'}$ is the excess vector bundle then
$$
\pi^*i_* ( a)= j_*(c_{n-m}(K)\pi '^*( a)).
$$
Moreover, we have the refined version
$$
\pi^*\p_i ( a)= \p_j(c_{n-m}(K)\pi '^*( a)).
$$
\end{coro}

\qed

\begin{obser}
Due to Corollary \ref{Coro inyectivo}, the \emph{key formula} and the fact that for a
regular closed immersion $i\colon Z\to X$ of codimension one we have $\p_i( a)= a
\cdot c_1^Z(L_Z)$ characterize the Gysin morphism and the refined Gysin morphism for
regular immersions. The general Gysin morphism is characterized by analogous
conditions (see upcoming Theorem \ref{Unicidad imagen directa} for the complete
statement).
\end{obser}

\subsection{Functoriality}

In order for the definition of the Gysin morphism to be of any use it has to be
functorial. In other words, if $Z\xrightarrow j Y \xrightarrow i  X$ are regular
immersions then the morphism $(i j)_*$ should be equal to $i_*j_*$. It is clear that
if the classes $ \bar \eta_Z ^X$ and $\bar \eta _Y ^X \bar \eta _Z ^Y \in \EE_Z (X)$
coincide this readily implies the functoriality.

\begin{teo}\label{Teo funtorialidad refinada}
If $Z \xrightarrow j Y \xrightarrow i X$ are two regular immersions then
\begin{equation}\label{Eq Funtorialidad}
\bar \eta _{Z}^X =\bar \eta _{Z}^Y \cdot \bar \eta _{Y}^X \in  \EE_Z(X).
\end{equation}
\end{teo}
\demo Let $n$ be the codimension of $j$ and $m$ that of $i$, that we may assume
constant. We split the proof into two parts:

\begin{lema}
With the preceding notations, if equation (\ref{Eq Funtorialidad}) holds for $m=1$
then it holds for any $m$.
\end{lema}

\demo The case $m=0$ is trivial so that assume $m>0$. Consider $B=B_YX$ the blow-up of
$Y$ in $X$ and denote $P=\PP (N_{Y/X})$, $P'=\PP (N_{Y /X} | _Z)$. We have the diagram
$$
\xymatrix{
 P'\ar[r]^{j'}\ar[d]_{p'} & P \ar[r] \ar[d]_p & B  \ar[d]^\pi \\
Z \ar[r] & Y \ar[r] & X }
$$
where both squares are cartesian, $P\to B$ is a regular closed immersion of
codimension 1 and $j'$ is of codimension $n$. Since the morphism
$$
\EE _Y ^{2(n+m),n+m}(X)\buildrel{\pi ^*}\over{\longrightarrow} \EE
_{P'}^{2(n+m),n+m}(B)
$$
is injective (\emph{cf.} Corollary \ref{Coro inyectivo}) it is enough to check in
$\EE _{P'}^{2(n+m),n+m}(B)$ the relation. Denote $\K^*$ the kernel of the epimorphism
$p^{'*} \I_{Z}/\I_{Z}^{ 2} \to \I_{P'} /\I_{P'} ^{2}$ and $K$ its associated vector
bundle. Using that the refined fundamental class are stable under base change
(Proposition \ref{Funtorialidad clase fund}), the formula from Lemma \ref{Lema
epimorfismo} and the equation (\ref{Eq Funtorialidad}) for $m=1$ we get
$$
\pi ^* \bar \eta _{Z}^X = \bar \eta _{P', \pi^* N_{Z/X}} ^{B}=c_{n-1} (K)\bar \eta
_{P'} ^{B}=c_{n-1} (K)\bar \eta _{P'} ^P \bar \eta _{P} ^B.
$$
Now, consider the commutative diagram of vector bundles on $P'$
$$
\begin{matrix}
& & & & 0 & & 0 & & \\
& & & & \uparrow & & \uparrow & &\\
& & 0 &\longrightarrow & K &\longrightarrow &
j'^* K' &\longrightarrow & 0\\
& & \uparrow & & \uparrow & &\uparrow & &\\
0&\longrightarrow &p'^* N_{Y/X}  &\longrightarrow & p'^* N_{Z/X}
&\longrightarrow & j'^* p^* N_{Y/X} & \longrightarrow & 0 \\
& & \uparrow & & \uparrow & &\uparrow & &\\
0&\longrightarrow &N_{P' / P} &\longrightarrow & N_{P'/B}
&\longrightarrow &  j'^* N_{P/B} & \longrightarrow& 0 \\
& & \uparrow & & \uparrow & &\uparrow & &\\
& & 0 & & 0& &0 & &
\end{matrix} $$
where $K'^*$ is the vector bundle associated to the kernel of $p^* \I_{Y}/\I_Y^2\to
\I_{P}/\I_{P}^2$. Taking into account that $j'^* K' \buildrel{\sim}\over{\to} K$ and
using Proposition \ref{Funtorialidad clase fund} and Lemma \ref{Lema epimorfismo}
once again we conclude

\begin{align*}
\pi ^* \bar\eta _{Z}^X &=c_{n-1} (K')\bar\eta _{P'} ^P \bar\eta _{P} ^{B}
=(c_{n-1}(K) \bar\eta _{P} ^{B})\bar\eta _{P'} ^P=
\bar\eta _{P,p^* \I_{Y}/\I_Y^2} ^{X_Y}\bar\eta _{P'} ^P \\
& =\pi ^* (\bar\eta _{Y} ^X) \bar\eta _{P'} ^P =\pi ^* \bar\eta _{Y} ^X \pi^*
\bar\eta _{Z} ^Y.
\end{align*}

\qed

\begin{lema}\label{Lema crucial funtorialidad}
With the preceding notations, the equation (\ref{Eq Funtorialidad}) is true for
$m=1$.
\end{lema}
\demo Denote $P=\PP (N_{Z/X})$, $P'= \PP (N_{Z/Y})$, $n=\codim _Y Z $ and $Y_Z$,
$X_Z$ for the blow-up of $Z$ in $Y$ and $X$ respectively. Consider the commutative
diagram
$$
\xymatrix{
P' \ar[r]^h \ar[d]_-v & Y_Z \ar[d] ^-w  & \\
P \ar[r]^-g \ar[d]& \pi ^{-1}(Y) \ar[r] \ar[d] & X_Z\ar[d]^-\pi \\
Z \ar[r]^j &Y \ar[r] & X }
$$
where every square is cartesian. Since $\pi^* \colon \EE_Z ^{2n+2,n+1}(X) \to \EE _P
^{2n+2,n+1}(X_Z)$ is injective (\emph{cf.} Corollary \ref{Coro inyectivo}) it is
enough to prove

\begin{equation}\label{Eq Funt subida}
\pi^* \bar \eta _Z  ^X = \pi ^* \bar \eta _Z ^Y \pi ^* \bar \eta _Y ^X
\end{equation}

\medskip

\noindent where $\pi^* \bar \eta _Z  ^X \in \EE _P ^{2n+2,n+1}(X_Z)$, $\pi ^* \bar
\eta _Z ^Y \in \EE _P ^{2n,n}(\pi^{-1}(Y))$ and $\pi ^* \bar \eta _Y ^X \in \EE
^{2,1}_{\pi ^{-1}(Y)}(X_Z)$.

We make the explicit computations of these two terms. This is the point where we
separate from \cite{Riou}, since we do not assume Chern classes to be additive as in
the case of \'{e}tale cohomology.

Let $\I_Y ^X$ be the sheaf of ideals of $Y$ in $X$ and $L_Y$ the line bundle
associated to its dual. To begin with,

\begin{align*}
\pi ^* \bar \eta _Y ^X &=c_1 ^{\pi ^{-1}Y}(\pi^*L_Y)=c_1 ^{\pi ^{-1}Y}(L_P \otimes
L_{Y_Z}) =c_1 ^{\pi ^{-1}Y}(\O_{X_Z}(-1)\otimes
L_{Y_Z})\\
&= c_1^{\pi ^{-1}Y}(\O_{X_Z}(-1))+c_1 ^{\pi ^{-1}Y}(L_{Y_Z})+c_1 ^{\pi
^{-1}Y}(\O_{X_Z}(-1))c_1 ^{\pi^{-1}Y}(L_{Y_Z})f\\
\end{align*}

\noindent where $f\in \EE^{**}(S)[[x,y]]$ is the series given by the formal group law
(\emph{cf.} \ref{serie f} and \ref{Propo serie f soporte}). Therefore, the right hand
side of equation (\ref{Eq Funt subida}) is the sum of the preceding three terms
multiplied by $\pi^* \bar \eta_Z ^Y$.

We compute each one of those three terms. From now on, we use the notation $u=c_1 ^P
(\O_{X_Z} (-1))\in \EE^{2,1} _P (X_Z)$. For the first term
$$
\pi^* \bar \eta _Z ^Y  \cdot c_1 ^{\pi^{-1}Y}(\O_{X_Z}(-1))= g^*\pi^* \bar\eta_Z ^Y
\cdot u= c_n(\pi^* N_{Z/Y})u = \mathrm{I}_1
$$
where the first equality is due to 6 of Proposition \ref{Propo propiedades producto}
applied to the map $g_\flat\colon \EE_P(X_Z)\to \EE_{\pi^{-1}Y}(X_Z)$ and the second
one to $g^*\pi^*=\pi^*j^*$ together with Proposition \ref{Propo propiedades Gysin}.
For the second term
\begin{multline*}
\pi^* \bar \eta _Z ^Y\cdot  c_1 ^{\pi ^{-1}Y}(L_{Y_Z}) = \pi^*\bar \eta_Z ^Y \cdot
w_\flat c_1 ^{Y_Z} (L_{Y_Z})=v_\flat(w^*\pi^* \bar \eta
_Z ^Y \cdot c_1 ^{Y_Z}(L_{Y_Z})) \\
= v_\flat \big((-1)^{n+1}[\sum _{i=0}^{n-1} (-1)^{i} c_i (N_{Z/Y})c_1 (\O _{P'}(-1))
^{n-1-i}]c_1 ^{P'}(\O _{Y_Z}(-1))c_1 ^{Y_Z}(L_{Y_Z})\big) \\
= v_\flat\big((-1)^{n+1}[\sum _{i=0}^{n-1} (-1)^{i} c_i (N_{Z/Y})c_1 ^{P'}(\O
_{Y_Z}(-1)\big)^{n-i}]c_1 ^{Y_Z}(L_{Y_Z}))\\
= (-1)^{n+1}[\sum _{i=0}^{n-1} (-1)^{i} c_i (N_{Z/Y})u^{n-i}c_1
(L_{Y_Z})]=\mathrm{I_2}
\end{multline*}
is due to Proposition \ref{Propo propiedades producto} for $w_\flat\colon
\EE_{Y_Z}(X_Z)\to \EE_{\pi^{-1}Y}(X_Z)$, Proposition \ref{Propo funtorialidad clase
Chern refi} and Corollary \ref{Coro c1 soporte}. For the third and last term
\begin{multline*}
\pi^* \bar \eta _Z ^Y  c_1 ^{\pi ^{-1}Y}(L_{Y_Z}) c_1 ^{\pi^{-1}Y} (\O _{X_Z}(-1))f\\
=(-1)^{n+1}[\sum _{i=0}^{n-1} (-1)^{i} c_i (N_{Z/Y})u^{n-i}]c_1 (L_{Y_Z})
c_1^{\pi^{-1}Y} (\O_{X_Z}(-1)) f\\
=(-1)^{n+1}[\sum_{i=0} ^{n-1}(-1)^{i} c_i (N_{Z/Y})u^{n-i+1}c_1(L_{Y_Z} )f]
=\mathrm{I}_3
\end{multline*}
where we use the preceding computation.

Consider the short exact sequence
$$
0\to \K^* \to \I_{Z}^X/\I_Z^{X2}\to \I_Z^Y /\I_Z^{Y2}\to 0
$$
With it, we compute the other side of the equation (\ref{Eq Funt subida}):
\begin{align*}
\pi ^* \bar \eta _Z ^X &= (-1)^{n}[\sum _{j=0} ^n (-1)^{j} c_j (N_{Z/X})c_1
(\O _P (-1))^{n-j})]u\\
&= (-1)^n[\sum _{j=0} ^n (-1)^{j} (c_j(N_{Z/Y})+c_{j-1}(N_{Z/Y})c_1(K))u^{n+1-j}]
\end{align*}
Note that $\I^P _{P'}=\K\otimes \O_P (-1)$ and therefore $\K=\I_{Y_Z}^{X_Z} \otimes
\O_P(1)$ so that $c_1(K)=c_1(L_{Y_Z}\otimes \O_P(-1))$.
\begin{align*}
=& \sum_{j=0}^n (-1)^{n+j}[ c_j(N_{Z/Y})+c_{j-1}(N _{Z/Y}) (c_1(L_{Y_Z})+c_1(\O_{P}
(-1))+c_1(L_{Y_Z})c_1(\O_{P}(-1))f)]u^{n+1-j}
\end{align*}
Therefore, $\pi ^* \bar \eta _Z ^X$ is the sum of three terms:
\begin{align*}
 &\sum_{j=0}^n (-1)^{n+j}
[c_j(N_{Z/Y})+c_{j-1}(N _{Z/Y}) c_1(\O_{P}(-1))]u^{n+1-j}=
\\
=&\sum_{j=0}^n (-1)^{n+j}c_j(N_{Z/Y})u^{n+1-j} + \sum_{i=0}^{n-1}
(-1)^{n+i+1}c_{i}(N _{Z/Y}) u^{n+1-i}=\\
=&c_n (N_{Z/Y})u = \mathrm{I}_1
\end{align*}
which is given by Proposition \ref{Propo propiedades Gysin} and the definition of the
Chern class,
\begin{align*}
&\sum_{j=0}^n (-1)^{n+j} c_{j-1}(N _{Z/Y}) (c_1(L_{Y_Z})))u^{n+1-j}=\\
=&(-1)^{n+1}[\sum _{i=0}^{n-1} (-1)^{i} c_i (N_{Z/Y})u^{n-i}c_1 (L_{Y_Z})] =
\mathrm{I}_2
\end{align*}
and finally
\begin{align*}
&\sum_{j=0}^n (-1)^{n+j} [c_{j-1}(N _{Z/Y}) c_1(L_{Y_Z})c_1(\O_{X_Z}(-1))f]u^{n+1-j}
=\\
=&(-1)^{n+1}[\sum_{i=0} ^{n-1}(-1)^{i} c_i (N_{Z/Y})u^{n+1-i}c_1(L_{Y_Z} )f] =
\mathrm{I}_3.
\end{align*}

\qed

\begin{ejem}\label{Ejem clase Thom}
Let $V\to X$ be a vector bundle of rank $n$ and $(\EE, c_1)$ be an oriented absolute
ring spectrum. The \textbf{Thom class} of $V$ is defined to be
$$
\mathbf{t}(V)\coloneqq \sum_{i=0}^n(-1)^i c_i(V)x^i\in \EE^{2n,n} (\bar V)
$$
where $x=c_1(\O _{\bar V}(-1))$ and $\bar V=\PP (V\oplus 1)$. It has being standard
in motivic homotopy theory since its beginning to define fundamental classes out of
Thom classes. More concretely, denote $s_0\colon X\to \bar V$ the zero section. Its
fundamental class was, by definition,  $\t(V)$. Therefore, the unicity of Gysin
morphisms in the context of regular schemes (\emph{cf.} \cite{Deglise}) proves that
for a regular scheme $X$ then
$$
\eta_X ^{\bar V}= \t(V).
$$
\end{ejem}

\begin{parra}\label{Propiedades clase de Thom}
Let us check that $\t(V)$ coincides with $\eta_X ^{\bar V}$ for arbitrary schemes. In
order to do so we recall some facts of the theory of Thom classes. For convenience of
the reader we recall the definitions. We define the \textbf{Thom space} of $V$ as
$$
\Th(V)=V/V-0\simeq \bar V /\PP(V).
$$
Its cohomology fits into a long exact sequence
$$
\ldots \to \EE ^{**} (\Th (V)) \xrightarrow{\pi^*} \EE^{**}(\bar V)\to \EE^{**}(V)
\to \ldots
$$
where, from Theorem \ref{Teo Hirsch Leray}, the third arrow is always a split
epimorphism. Since $\t (V)$ is zero in $\EE(\PP(V))$, we call the \textbf{refined
Thom class} to the unique element
$$
\bar \t (V)\in \EE(\Th (V))\simeq \EE_X (\bar V)=\EE_X(V)
$$
such that $\pi^*(\bar \t (V) )=\t (V)$. Clearly, proving that $\bar \t (V)$ coincides
with $\bar \eta_X ^{\bar V}$ is equivalent to proving that $\bar \t (V)$ coincides
with $\bar \eta_X ^{\bar V}$.

One last technical recall (\cf \cite{Deglise} for example): if $0\to V'\to V\to
V''\to 0$ is exact, the refined Thom classes satisfy
$$
\bar \t(V)=\bar \t(V')\bar\t (V_{V'})\in \EE _X(V).
$$
Here $\t (V_{V'})$ denotes $V$ considered as a bundle over $V'$ \footnote {The scheme
which parametrises the sections of $V\to V''$ is a torsor (of group $\Hom (V'',V')$)
and  the pullback of the short exact sequence is naturally split there.} and the
product is that of Definition \ref{Defi producto refinado}.
\end{parra}

\begin{propo}\label{Clase de Thom}
Let $V\to X$ be a vector bundle and denote $ \bar V$ its projective completion. Then
$$
\eta _X^{\bar V}=\t (V).
$$
\end{propo}
\demo It is clear that this formula is equivalent to its refined counterpart, $\bar
\eta _X^{\bar V}=\bar \t (V)$. Due to Theorem \ref{Teo funtorialidad refinada} and
the previous remark it is enough to prove it for the case of a line bundle $V=L$.

In this case $\t (L)= c_1(L)-c_1(\O_{\bar L}(-1))$ and $\eta_X^{\bar L}=c_1(\I^*)$,
where $\I$ stands for the sheaf of ideals of the zero section in $\bar L$.  This
sheaf may be computed explicitly: the composition $\O_{\bar L}(-1) \to L\oplus \O \to
L$ of the canonical morphism and the projection is an isomorphism out of the zero
section, which induces  $ L^*\otimes \O_{\bar L}(-1)\simeq \I \to \O$. Consider the
canonical short exact sequence
$$
0\to \O_{\bar L}(-1) \to L\oplus \O \to Q\to 0
$$
where $Q$ is the canonical quotient bundle. Hence $L= \bigwedge^2 (L\oplus
\O)=Q\otimes \O_{\bar L}(-1) $ so that $Q= L\otimes \O_{\bar
L}(1)=\I^*$.\footnote{Note that we already computed a more general formula in the
proof Lemma \ref{Lema crucial funtorialidad} from which this equality follows.}

To conclude note that $c_1(L\oplus \O)=c_1(L)$. Also by the additivity of Chern
classes $ c_1(L\oplus \O)=c_1(Q)+c_1(\O(-1))$.\footnote{I am thankful to J. Riou who
pointed me out a mistake in this computation on an earlier version of this text.}

\qed

Let $\EE$ be an absolute ring spectrum and $\MM$ be an absolute $\EE$-module. Note
that $\MM$ does not have a unit. As a consequence there are no fundamental nor Chern
classes in the $\MM$-cohomology. However, $\MM(X)$ is an $\EE(X)$-module and
therefore we can still multiply classes in the $\MM$-cohomology by fundamental
classes and Chern classes of the $\EE$-cohomology. This suffices to define the Gysin
morphism.

\begin{defi}
Let $i\colon Z\to X$ be a regular immersion. We define the \textbf{Gysin morphism}
$i_*$ and the \textbf{refined Gysin morphism} $\p_i$ in the $\MM$-cohomology to be
$$
\begin{matrix}
i_*\colon  & \MM^{p,q}(Z) & \longrightarrow & \MM^{p+2n,q+n}(X) \\
&m & \mapsto & m\cdot \eta_Z^X
\end{matrix}
\quad , \quad
\begin{matrix}
\p_i\colon  & \MM^{p,q}(Z) & \longrightarrow & \MM^{p+2n,q+n}_Z(X) \\
&m & \mapsto & m\cdot \bar \eta_Z^X.
\end{matrix}
$$
\end{defi}

We readily deduce from the case of ring spectra the following properties for modules:

\begin{teo}\label{Teo propiedades Gysin modulos}
Let $\EE$ be an oriented absolute ring spectrum, $\MM$ an absolute $\EE$-module and
$i\colon Z \to X$ a regular immersion.
\begin{enumerate}

\item
Functoriality: Let $j \colon Y\to Z$ be a regular immersion, then $(ij)_* =i_* j_*$.

\item
Projection formula: The Gysin morphism is $\EE(X)$-linear. In other words,
$$
a\cdot i_* (m)= i_* (i^* (a) \cdot m) \qquad \forall \ a \in \EE (X) \ , \ m \in \MM
(Z).
$$
Note that an analogous formula also holds for $n \in \MM(X)$ and $b\in \EE(Z)$.

\item
Denote $n=\codim _X Z$, we have
$$
i^*  i_*(m) =c_n(N_{Z/X})\cdot m \quad \forall \ m\in \MM^{2n,n} (Z).
$$

\item
Let $r\colon X\to  Z $ be a retraction of $i$, then the Gysin morphism $i_*$ is $\MM
(Z)$-linear (with $r^*$). In other words,
$$
i_*(m)=r^*(m)\cdot \eta_Z ^X \quad \forall \ m \in \MM(Z).
$$

\item
Excess intersection formula: Consider a cartesian square
$$
\xymatrix{ P \ar[r]^j \ar[d]_{\pi '} &X'\ar[d]^\pi \\
Z\ar[r]^i & X}
$$
where both $i$ and $j$ are regular immersions of codimension $n$ and $m$
respectively. Denote $K=\pi'^* N_{Z/X}/N_{P/X'}$ the excess vector bundle, then
\begin{align*}
\pi^*i_* (m)&= j_*(c_{n-m}(K)\cdot \pi '^*(m)) \qquad \mbox{and} \\
\pi^*\p_i (m)& = \p_j(c_{n-m}(K)\cdot \pi '^*(m)) \qquad \forall \ m \in \MM (Z).
\end{align*}

\end{enumerate}

\end{teo}

\qed

\subsection{The projective lci case}

The main reference we have used for this section is $\cite[\S 5]{Deglise2}$ where
D\'{e}glise thoroughly studied the Gysin morphism. However, the reference works on the
smooth case and on a general category of premotives satisfying certain axioms, which
do not hold for $\SH$. Nevertheless, the arguments still hold \emph{mutatis mutandis}
in our context.

We first observe that the fundamental class of the diagonal
$\eta_{\Delta_n}=\eta_{\Delta_n}^{\PP^n_X\times \PP^n_X}\in \EE(\PP^n_X)\otimes \EE
(\PP^n_X)$ defines a nondegenerate bilinear pairing. The duality it defines (similar
to Euclidian spaces) allows us to define the Gysin morphism for the projection
$p_{{\scriptscriptstyle X}}\colon \PP^n_X \to X$ of a projective space onto its base
as the dual of the inverse image $p_{{\scriptscriptstyle X}}^*\colon \EE(X)\to
\EE(\PP^n_X)$ (\cf Definition \ref{Gysin proyeccion}). We then prove that the Gysin
morphism for projective lci morphism, without smoothness assumptions, have all usual
properties. We deduce from the case of ring spectra a Gysin morphism for modules.

\begin{parra}\label{dualidad}
Denote $\EE(\PP^n_X)^\vee= \Hom_{\EE(X)\mod}$ $(\EE (\PP^n_X),\EE(X))$ the dual of
$\EE (\PP^n_X)$ as $\EE(X)$-module, $p_{{\scriptscriptstyle X}}\colon \PP^n_X\to X$
the canonical projection and $$(p_{{\scriptscriptstyle X}}^*)^\vee\colon
\EE(\PP^n_X)^\vee\to \EE(X)^\vee\simeq \EE(X)$$ the transpose of the inverse image.
Recall that from the projective bundle theorem we have $ \EE (\PP^n_X\times_X
\PP^n_X)\simeq \EE(\PP^n_X)\otimes_{\EE(X)} \EE (\PP^n_X)$. Denote $\Delta_n \colon
\PP^n_X \to \PP^n_X\times \PP^n_X$ the diagonal embedding. The fundamental class of
the diagonal $\eta_{\Delta_n}\in \EE(\PP^n_X)\otimes \EE (\PP^n_X)$ defines a
symmetric $\EE(X)$-bilinear pairing on $\EE(\PP^n_X)^\vee$, and therefore a polarity
$$
\Phi \colon\EE(\PP^n_X)^\vee \to \EE(\PP^n_X) \mbox{ , } \omega \mapsto \
(\omega\otimes 1 )(\eta_{\Delta_n})=\eta_{\Delta_n}(\omega,\raya).
$$
\end{parra}

\begin{propo}
The polarity $\Phi$ defined above is an isomorphism.
\end{propo}
\demo The matrix that defines $\Phi$ is the same as $\eta_{\Delta_n}$. One can check
that, as in the smooth case, we have
\begin{equation}\label{Mn}
\eta_{\Delta_n}=\sum_{r,s=0}^n a_{r,s}x_n^r\otimes x^s_n \quad \mbox{where} \quad
(a_{r,s})= \left( \raisebox{0.5\depth}{ \xymatrix@=1ex{ 0\ar@{.}[rr]\ar@{.}[dd]
 & & 0\ar@{.}[lldd]
 & 1\ar@{-}[lllddd] \\
 & & & \bullet\ar@{.}[lldd]\ar@{.}[dd] \\
0 & & & \\
1 & \bullet\ar@{.}[rr] & & \bullet } } \right)
\end{equation}
where $x_n=c_1(\O_{\PP^n_X}(-1))$ (\emph{cf}. \cite[2.2]{Navarro Grothendieck}).

\qed

\begin{defi}\label{Gysin proyeccion}
Let $\EE$ be an oriented absolute ring spectrum and denote $p_{{\scriptscriptstyle
X}}\colon \PP^n_X\to X$ the natural projection, we define the \textbf{direct image}
of $p_{{\scriptscriptstyle X}}$ to be
$$
p_{{\scriptscriptstyle X}*}\colon  \EE(\PP^n_X)\buildrel{\Phi^{-1}}\over{\simeq
}\EE(\PP^n_X)^\vee \xrightarrow{(p^*)^\vee}\EE(X)^\vee=\EE(X).
$$
\end{defi}

\begin{obser}
Denote $p\colon \PP^n\to S$, since the square
$$
\xymatrix{\PP^n_X\ar[d]\ar[r]&\PP^n_X\times \PP^n_X\ar[d]^\pi\\
\PP^n\ar[r]&\PP^n\times \PP^n}
$$
is transversal, the excess intersection formula gives
$\eta_{\Delta_{nX}}=\pi^*\eta_{\Delta_n}$ and we have
$$
p_{{\scriptscriptstyle X}*}= p_*\otimes 1_X.
$$
\end{obser}

We specify a computation from linear algebra.

\begin{lema}
Let $A$ be a ring and $M$ be and finitely generated free $A$-module. Let $f\colon
A\to M$ be a linear map and $\eta\in M\otimes_A M$ a symmetric bilinear pairing on
$M^\vee$. If $\eta$ is nonsingular (in the sense that its associated polarity $\Phi
\colon M^\vee \to M \, , \, \alpha \mapsto (\alpha \otimes
1_M)(\eta)=\eta(\alpha,\raya)$ is an isomorphism) then the map
$$
\omega \colon M\xrightarrow {\Phi^{-1}}M^\vee\xrightarrow{f^\vee}A^\vee\xrightarrow
\sim A
$$
satisfies
$$
(\omega \otimes 1_{M})(\eta)=f(1).
$$
\end{lema}
\qed

\begin{obser}\label{Determinar p_*}
The previous lemma applied to our case $\eta=\eta_{\Delta_n}$ and
$f=p^*_{\scriptscriptstyle X}\colon \EE (\PP^n_X)\to \EE(X)$ concludes
$$
\Phi(p_{{\scriptscriptstyle X}*})=p_{\scriptscriptstyle X}^*(1)=1,
$$
which totally determines $p_{{\scriptscriptstyle X}*}$.
\end{obser}

\begin{lema}
Let $\Delta_n \colon \PP^n\to \PP^n\times \PP^n=\PP^n_{\PP^n}$ be the diagonal
embedding. Then $p_{2*}\Delta_*=1_{\PP_n}$.
\end{lema}
\demo Recall that $\{1, x=c_1(\O_{\PP^n}(-1))$, $\ldots$ , $x^n\}$ is a basis of $\EE
(\PP^n)$. From last lemma applied to $f=p^*$, $\omega=p_*$ and $\eta=\Delta_*(1)$ we
have
$$
1=p^*(1)=(p_*\otimes 1)(\eta_{\Delta_n})=p_{2*}(\Delta_*(1)).
$$
For the rest of elements of the basis we conclude by observing that
$$
p_{2*}(\Delta_*(x^j))=p_{2*}(\Delta_*(\Delta^*(1\otimes x^j))=p_{2*}(1\otimes
x^j)\cdot p_{2*} (\Delta_*(1))=x^j
$$
since $p_{2*}=p_*\otimes 1_{\PP^n}$.

\qed

\begin{coro}
Let $s\colon X\to \PP^n_X$ be a section of $p_{\scriptscriptstyle X}\colon \PP^n_X\to
X$. Then $p_{{\scriptscriptstyle X}*}s_*=1_X.$
\end{coro}
\demo Note that preceding lemma also holds for a projective space over a general base
$\PP^n_X$. Consider the cartesian squares
$$
\xymatrix{X \ar[r]^-s \ar[d]_{s}& \PP^n\times X \ar[r]^-{p_{X}} \ar[d]^-{1_{\PP^n}
\times s} & X \ar[d]^{s}\\
\PP^n_X\ar[r]^-\Delta & \PP^n\times \PP^n \times X\ar[r]^-{p_{ \PP^n_X}} & \PP^n_X.}
$$
We may apply Corollary \ref{Coro formula exceso} to the left square so that
$s^*s_*=\Delta_*(1_{\PP^n_X}\times s)^*$. From Definition \ref{Gysin proyeccion} we
also have that $p_{{\scriptscriptstyle X}*}(1_{\PP^n_X}\times s)^*=p_{\PP^n_X*}$.
Together with the previous lemma we deduce $s^*=p_*s_*s^*$. Since $s^*$ is surjective
we conclude that $p_*s_*=1_X$.

\qed

\begin{lema}\label{Preparacion defi Gysin lci}
Let $i\colon Z\to  X$ be a regular immersion and consider the cartesian diagram
$$
\xymatrix{\PP^n _Z \ar[r]^k \ar[d]_{p_Z} &\PP^n_X \ar[d]^{p_X} \\
Z \ar[r]^i & X}
$$
Then $k_*=( 1_{\PP^n}\times i )_* = 1_{\PP^n}\otimes i_* $ and $
p_{{\scriptscriptstyle X}*}k_*=i_*p_{{\scriptscriptstyle Z}*}$.

\end{lema}
\demo For the first claim, applying the excess intersection formula and the
projective bundle theorem (Corollary \ref{Coro formula exceso}, Theorem \ref{Teo
Hirsch Leray}) we get that the diagram
$$
\vcenter{\xymatrix{ \EE (\PP^n)\otimes_{\EE(S)} \EE(Z) \ar[r]^{k_*} &
\EE (\PP^n)\otimes_{\EE (S)} \EE(X)\\
\EE(Z) \ar[r]^-{i_*} \ar[u]^{ p_Z^*} &\EE (X) \ar[u]_{p_X^*}}}
$$
commutes so that $k_*(1\otimes  a)=1\otimes i_*( a)$ for $ a \in \EE(Z)$. Applying
the projection formula from Proposition \ref{Propo propiedades Gysin} we get that
$k_*( b \otimes 1  )=( b \otimes 1)\cdot k_*(1\otimes 1)= b\otimes \eta_Z ^X$ for $ b
\in \EE(\PP^n)$ so we conclude $k_*=(1_{\PP^n}\times i  )_* = 1_{\PP^n}\otimes i_* $.
From here and the previous definition the formula $p_{{\scriptscriptstyle
X}*}k_*=i_*p_{{\scriptscriptstyle Z}*}$ follows.

\qed

\begin{teo}\label{Teo lci no depende factorizacion}
Consider a commutative diagram
$$
\xymatrix{
Y \ar[r]^k \ar[d]_i &\PP^n _X \ar[d]^p\\
\PP^m _X \ar[r]^q & X }
$$
where $i$ and $k$ are regular immersions of codimension $r$ and $s$ respectively and
$p$ and $q$ are the natural projections. Then, $p_* k_*  = q_*i_* $.
\end{teo}
\demo Consider the following commutative diagram:
$$
\xymatrix@C=22pt@R=8pt{
&& & {}\PP^n_X\ar^p[rd] & \\
Y\ar@/^5pt/^k[rrru]\ar@/_5pt/_i[rrrd]\ar|v[rr]
 && {}\PP^n_X \times_X {}\PP^m_X\ar_/8pt/{q'}[ru]\ar^/8pt/{p'}[rd]
 & & X. \\
&& & {}\PP^m_X\ar_q[ru] & }
$$
Since it is clear that $p_*q'_*=q_*p'_*$ it is enough to prove that $p'_*v_*=i_*$.
For that case, denote $T=\PP^m_X$, $v=\rho \times i$ where $\rho \colon Y\to \PP^n$
and consider
$$
\xymatrix@R=5pt@C=7pt{
Y\ar@{=}[rrddd]\ar^/4pt/s[rrd]\ar@/^8pt/^v[rrrrd] && & \\
&& {}\PP^n_Y\ar_l[rr]\ar^\pi[dd] && {}\PP^n_T\ar^{p'}[dd] \\
&& && \\
&& Y\ar^i[rr] && T }
$$
where $i$, $s=\rho\times 1_Y$, $l$ and $v$ are regular immersions. By the
functoriality of Theorem \ref{Teo funtorialidad refinada} we have $v_*=l_*s_*$ and by
the previous Lemma \ref{Preparacion defi Gysin lci} we also have $\pi_*s_*=1_Y$ and
$p'_*l_*=i_*\pi_*$. Considering all together we conclude
$$
p'_*v_*=p'_*l_*s_*=i_*\pi_*s_*=i_*.
$$

\qed

\begin{defi}\label{Defi proyectivo lci}
We define an $X$-scheme $Y\to X$ to be a \textbf{local complete intersection}
(\textbf{lci}) if it locally admits a factorization by a regular immersion into
$\AA^n_X$ (\cite[VIII 1.1]{SGA6}). Given our conventions, a projective lci morphism
$f\colon Y\to X$ admits a factorization of the form $Y\xrightarrow i\PP^n _X
\xrightarrow p X$ where $i$ is a regular closed immersion (\cf \cite[VIII 1.2]{SGA6})
and $p$ is the canonical projection.

Let $\EE$ be an absolute ring spectrum, $f\colon Y\to  X$ be a projective lci
morphism and $Y\xrightarrow i \PP^n _X \xrightarrow p Y$ be a factorization. We
define the \textbf{direct image} of $f$ as $f_*\coloneqq p_*i_*$ (by Theorem \ref{Teo
lci no depende factorizacion} it does not depend on the choice of factorization).

Let $\MM$ be an absolute $\EE$-module and $p_{\scriptscriptstyle X}\colon \PP^n_X\to
X$ be the natural projection. We define the direct image $p_{{\scriptscriptstyle
X}*}$ in the $\MM$-cohomology as the morphism
$$
\xymatrix@R=6pt{\MM(\PP^n_X)\ar@{-->}[rr]^-{p_{X*}}\ar@{=}[d]&&\MM(X)\ar@{=}[d]\\
\EE(\PP^n)\otimes_{\EE(S)}\MM(X)\ar[rr]^-{p_*\otimes 1}&& \EE(S)\otimes_{\EE(S)}
\MM(X).}
$$
Let $f\colon Y\to X$ be a projective lci morphism and $f=pi$ be a factorization
$Y\xrightarrow i \PP^n _X \xrightarrow p X$ where $i$ is a regular immersion and $p$
is the natural projection. We define the direct image of $f$ as $f_*\coloneqq
p_*i_*$.
\end{defi}

Finally, let us remark the main properties the direct image in this context.

\begin{teo}
Let $f\colon Y\to  X$ be a projective lci morphism, $\EE$ be an absolute oriented
ring spectrum and $\MM$ be an absolute $\EE$-module:

\begin{enumerate}
\item
Functoriality: If $g\colon Z\to Y$ is another projective lci morphism then
$$
(fg)_*=f_*g_*\colon \MM(Z)\longrightarrow \MM(X).
$$

\item
Excess intersection formula: Consider a cartesian square
$$
\xymatrix{Y'\ar[r]^-g \ar[d]_-q & X' \ar[d]^-p\\
Y\ar[r]^-f & X}
$$
where $f$ and $g$ are projective lci morphisms of codimension $r$ and $s$
respectively. Choose a factorization $Y\xrightarrow i\PP^n _{X} \xrightarrow p X$ of
$f$ and denote $K=q^*N_{Y/\PP^n_X}/N_{Y'/\PP^n_{X'}}$.\footnote{Recall that the
definition of $K$ in the previous proposition does not depend on the choice of
factorization (\emph{cf.} \cite[6.6]{Fulton}).} Then
$$
p^*f_*(m)=g_*(c_{r-s}(K)\cdot q^*(m)) \quad \forall \ m \ \in \ \MM(Y).
$$

\item
Projection formula: $f_*$ is a morphism of $\EE(X)$-modules. That is to say,
$$
f_* (f^*(a)\cdot m)=a\cdot f_* (m) \quad  \forall \ a \in \EE (X) \ , \ m \in \MM(Y).
$$

\end{enumerate}
\end{teo}
\demo We prove the case of ring spectra and the case of modules follows
straightforwardly.

For the first point consider factorizations $Z\xrightarrow j \PP^m_X \xrightarrow q
X$ of $f\circ g$ and $Y\xrightarrow i \PP^n_X \xrightarrow p X$ of $f$. We can
compute explicitly the base change
$$
\xymatrix{
Y'\ar[r]^{i'}\ar[d]_{\pi'} & \PP^m_X\times \PP^n_X\ar[d]^{\pi}\\
Y \ar[r]^i&\PP^n_X}
$$
as $Y'= (\PP^m_X\times _X \PP^n_X)\times_{\PP^n_X}Y=\PP^m_Y$. If we denote $j=v\times
(fg) \colon Z \to \PP^n\times X$ and consider $k=v\times g\colon Z\to \PP^n \times Y$
it fits into a commutative diagram
$$
\xymatrix@C=16pt@R=6pt{
&& {}\PP^m_X\ar^q@/^18pt/[rrddd] && \\
& & {}\PP^m_X \times {}\PP^n_X\ar^/5pt/{\pi}[rd]\ar_{p'}[u] & & \\
& {}\PP^m_Y\ar^{\pi'}[rd]\ar^/-6pt/{i'}[ru] & & {}\PP^n_X\ar|p[rd] & \\
Z\ar|g[rr]\ar|k[ru]\ar^j@/^18pt/[rruuu] && Y\ar|f[rr]\ar^i[ru] && X. }
$$
The preceding lemmas allow to conclude the functoriality since
$$
f_*g_*=p_*i_*\pi'_*k_*=p_*\pi_*i'_*k_*=q_*p'_*i'_*k_*=q_*j_*=(fg)_*.
$$

For the excess intersection formula recall the factorization $Y\xrightarrow i\PP^n
_{X} \xrightarrow p X$. Changing base on the regular immersion we get a cartesian
diagram
$$
\xymatrix{Y'\ar[r]^-j \ar[d]_-q & \PP^n_{X'}
 \ar[d]^-{p'}\\
Y \ar[r]^-i & \PP ^{n} _{X}}
$$
which has the same excess bundle $K$ and where $j$ and $i$ are regular immersions.
Hence by Corollary \ref{Coro formula exceso} for any $ a \in \EE(Y)$ the relation
$$
\pi^*i_*( a)=j_*(c_{n-m}(K)\cdot q^*( a))
$$
holds. Consider the diagram
$$
\xymatrix{\PP^n_{X'}\ar[r]^{\pi'}\ar[d]_{p'} & X' \ar[d]^-p\\
\PP^n_X\ar[r]^\pi & X}
$$
then $p^*\pi_*=\pi'_*p^{'*} $ and we conclude.

For the projection formula consider the commutative diagram
$$
\xymatrix{Y\ar[r]^f \ar[d]_{\gamma_f} & X\ar[d]^\Delta\\
Y\times X\ar[r]^{f\times 1_X}& X\times X}
$$
where $\gamma_f$ denotes the graphic of $f$ and $\Delta$ denotes the diagonal. Since
$\Delta$ is transversal to $f\times 1_X$ we may apply the excess intersection
formula: for any $a \in \EE(X)$ and $b \in \EE(Y)$ we have that
$$
\Delta^*(f\times 1_X)_*(b \times a)=\Delta^* ((f_*b )\times a)=f_*(b) \cdot a
$$
equals
$$
f_*\gamma_f^*(b\times a)=f_*(b\cdot f^*(a)).
$$

\qed

\medskip

As a result of the construction we can characterize direct images. The result
develops \cite[2.5]{Panin2} for smooth $k$-schemes and \cite[3.3.1]{Deglise} for
regular schemes.

\begin{teo}\label{Unicidad imagen directa}
Let $(\EE, c_1)$ be an oriented absolute ring spectrum, there exists a unique way of
assign for any projective lci morphisms $f\colon Y\to X$ a group morphism $f_*\colon
\EE (Y) \to \EE(X)$ satisfying the following properties:
\begin{enumerate}

\item
Functoriality: $(fg)_*=f_*g_*$.

\item
Normalization: For regular immersions $i\colon Y\to X$ of codimension one they
satisfy $i_*(a)=i_\flat (a \cdot c_1^Y(L_Y))$.

\item
Key formula: If $i \colon Y\to X$ is a regular immersion of codimension $n$,
$\pi^*\colon B_YX\to X$ is the blowing-up of $Y$ in $X$ with exceptional divisor
$j\colon \PP(N_{Y/X})\to B_YX$, we have $ \pi^* i_* (a)= j_*(c_{n-1} (K)\cdot \pi
'^*(a)). $

\item
Projection formula: They are $\EE(X)$-linear, i.e., $ f_*(f^* (a) \cdot b )= a \cdot
f_*(b)$ for $a \in \EE(X)$ and $b \in \EE(Y)$.

\end{enumerate}

When considering regular immersions, properties 2 and 3 characterize them.
\end{teo}
\demo The functoriality property reduces the proof to the case of regular immersions
and the projection of a projective space onto its base. The case of closed immersions
follows directly from properties 2, 3 and the long exact sequence of the blow-up (\cf
Corollary \ref{Coro inyectivo}).

For the projection of a projective space we apply the excess formula to the
commutative square
$$
\xymatrix{\PP^n_X\ar[r]^{p_X}\ar[d]_{\pi_1}&X\ar[d]\\
\PP^n\ar[r]^p&S.}
$$
Together with the projective bundle theorem \ref{Teo Hirsch Leray} and the projection
formula we obtain that $p_{{\scriptscriptstyle X}*}=p_*\otimes 1_X$. The
functoriality property implies the commutativity of the triangle
$$
\xymatrix{\EE(\PP^n)\ar[r]^-{\Delta_*}\ar[rd]_{1_{\PP^n_X}}&\EE(\PP^n)\otimes
\EE(\PP^n)\ar[d]^{p_*}\\
&\EE(\PP^n).}
$$
The argument from Remark \ref{Determinar p_*} shows that $p_*\in \EE(\PP^n)^\vee$ is
uniquely determined by its image through the polarity, which is
$\Phi(p_*)=(p_*\otimes 1 )(\eta_\Delta)=1$.

\qed

\begin{obser}
By the previous theorem the direct image constructed here coincides with classic
constructions for cohomologies represented by spectra of \ref{Ejem espectros en
anillo} (for instance \emph{cf.} \cite[3.3.4]{Deglise}). Let us explicitly mention
the case of higher $K$-theory: all properties apart from normalization are proved for
Quillen's $K$-theory in \cite{Quillen} and normalization is proved in
\cite{Thomason}.
\end{obser}

\section{Riemann-Roch theorem}

We devote this section to prove the motivic Riemann-Roch theorem in the context of
the algebraic stable homotopy category. We deduce a Riemann-Roch theorem for modules.
The classic Chern character is replaced in the general setting by a morphism of
oriented spectra $\varphi \colon  (\EE, c_{1}) \to (\FF , \bar c_{1})$. For clarity
in the exposition, overlined morphisms and elements will refer to the
$\FF$-cohomology.

\goodbreak

\begin{teo}\label{RR Facil}
Let $\varphi\colon  (\EE , c _1) \to (\FF , \bar c _1)$ be a morphism of oriented
absolute ring spectra such that $\varphi_{\PP^{\infty}} (c_{1})=\bar c_{1}$ and let
$f\colon Y\to X$ be a projective lci morphism, then the diagram
$$
\xymatrix{
\EE (Y) \ar[r]^{f_*}\ar[d]_{\varphi_Y} & \EE (X)\ar[d]^{\varphi_X}\\
\FF (Y) \ar[r]^{\bar f_*} & \FF (X)}
$$
commutes. In other words, for $a \in \EE(Y)$
$$
\varphi_X(f_*(a))= \bar f_* (\varphi_Y (a)).
$$
\end{teo}

\demo Following the standard approach, it is enough to check the theorem for a
regular immersion $i\colon X\to \PP^n_Y$ and the projection $p\colon \PP^n_Y \to Y$.

\medskip

\begin{lema}[Regular immersions]
Theorem \ref{RR Facil} holds for a regular immersion.
\end{lema}
\demo Since $\varphi$ preserves the orientation it also preserves Chern classes of
vector bundles. It follows that $\varphi (\eta_Z^X)=\bar \eta_Z^X$ since fundamental
classes were defined it terms of Chern classes. We conclude that $ \varphi_X i_*=
\bar i_* \varphi_Z $. \footnote{I am grateful to J. Riou for an observation that
simplified this proof as it is now.}

\qed

\begin{lema}[Projection]
Theorem \ref{RR Facil} holds for the canonical projection $p\colon \PP^n_X\to X$.
\end{lema}
\demo Applying Theorem \ref{RR Facil} to the diagonal embedding $\Delta_n\colon
\PP^n_X\to \PP^n_X \times \PP^n_X$ we obtain that $\varphi_{\PP^n_X \times \PP^n_X}$
preserves the fundamental class of the diagonal. Recall from Definition \ref{Gysin
proyeccion} that $p_*=\Phi^{-1}(p^*)^\vee$ where $\Phi\colon \EE(\PP^n_X)\xrightarrow
\sim \EE(\PP^n_X)^\vee$ is the polarity defined by the fundamental class of the
diagonal (\cf Paragraph \ref{dualidad}). Since $\varphi$ commute with inverse images
the diagram
$$
\xymatrix{ \EE(\PP^n_X) \ar[r]^{\Phi^{-1}}\ar[d]_{\varphi_{\PP^n_X}}
&\EE(\PP^n_X)^\vee \ar[d]^{\varphi_ {\PP^n_X}^\vee} \ar[r]^-{(p^*)^\vee} &\EE(X)^\vee
\simeq\EE(X)\phantom{aaaaaa}\ar[d]^-{\varphi_X}\\
\FF(\PP^n_X) \ar[r]^-{\bar \Phi^{-1}}&\FF(\PP^n_X)^\vee\ar[r]^-{(\bar
p^*)^\vee}&\FF(X)^\vee\simeq \FF(X) \phantom{aaaaaa}}
$$
is made of commutative squares.

\qed

\medskip

\begin{lema}[Change of direct image]\label{Nuevas imagenes directas}
Let $(\EE, c_1)$ be an oriented absolute ring spectrum. Let
$c_1^{\mathrm{new}}=G(c_1)\cdot c_1$ be a new orientation (cf. Proposition \ref{serie
c1}) and denote $G^{-1}_\times$ the multiplicative extension of
$G^{-1}\in\EE(S)[[t]]$ (cf. Paragraph \ref{Tod}). Let $f\colon Y\to X$ be a
projective lci morphism and denote $T_f=i^*T_{\PP^n_Y}-N_i\in K_0(Y)$ the virtual
tangent bundle of $f=p\circ i$, then
$$
f_*^{\mathrm{new}}(a)=f_*\bigl(G^{-1}_\times (T_f)\cdot a\bigr) \quad \forall \ a \
\in \EE(Y).
$$
\end{lema}
\demo We have to check that the right hand side of the equality satisfies the
conditions of Theorem \ref{Unicidad imagen directa}. In the case of the canonical
projection $p\colon \PP^n_ X \to X$ we only have to check the projection formula,
which is immediate. For the case of regular immersions $i\colon Y\to X$ we have to
check that the family of morphism $i_*^{\mathrm{new}}\colon \EE^{*,*}(Y) \to
\EE^{2d+*,d+*}(X)$, $a \mapsto i_*\bigl(G^{-1}_\times(-N_{Y/X})\cdot a\bigr)$ satisfy
the properties 2 and 3 of Theorem \ref{Unicidad imagen directa}. The first property
follows directly from the definition of the multiplicative extension. For the key
formula consider Corollary \ref{Coro inyectivo} notations and recall that we denote
the canonical quotient bundle by $K=\pi'^* N_{Y/X}/\O_P(-1)$. We then have
\begin{align*}
\pi' i_*^{\mathrm{new}}(a)&=\pi^*i_*\bigl(G^{-1}_\times(-N_{Y/X})\cdot a\bigr)=j_*
\bigl(c_{d-1} (K)\cdot
G^{-1}_\times (-K-\O_P(.-1)) \cdot \pi^* (a )\bigr) \\
&= j_*^{\mathrm{new}} \bigl(c_{d-1}^{\mathrm{new}}(K)\cdot \pi^*(a)\bigr)
\end{align*}
where the last equality comes from Corollary \ref{Tod cn}.

\qed

\begin{teo}[Motivic Riemann-Roch] \label{Motivic Riemann Roch}
Let $\varphi\colon  (\EE , c_1) \to (\FF , \bar c_1)$ be a morphism of oriented
absolute spectra. Denote $G\in \FF[S][[t]]$ the series such that $\varphi
(c_{1})=G(\bar c_1)\cdot \bar c_1\in \FF (\PP_S ^\infty)$ and $G^{-1}_\times$ be the
multiplicative extension of $G^{-1}$. Let $f\colon Y\to  X$ be a projective lci where
$f=p\circ i$ and denote $T_f= i^{*}T_p-N_i\in K_0(Y)$ the virtual tangent bundle.
Then the diagram
$$
\xymatrix{
\EE(Y)\ar[r]^{f_*} \ar[d]_{G^{-1}_\times(T_f)\varphi_Y}& \EE(X)\ar[d]^{\varphi_X}\\
\FF (Y) \ar[r]^{\bar f_*} & \FF (X) }
$$
commutes. In other words, for $a \in \EE(Y)$ we have
$$
\varphi_X(f_*(a))=\bar f_*\bigl(G^{-1}_\times(T_f)\cdot \varphi_Y(a)\bigr).
$$
\end{teo}
\demo We define $\bar c_1 ^{\mathrm{new}}=\varphi_{\PP^\infty}( c_{1})=G(\bar
c_1)\cdot \bar c_1$ that gives a direct image $\bar f_* ^{\mathrm{new}}$ satisfying
$\varphi_X (f_*(c_1))=\bar f_*^{\mathrm{new}}(\varphi_Y(\bar c_1))$ due to Theorem
\ref{RR Facil}. We conclude recalling Lemma \ref{Nuevas imagenes directas}.

\qed

\begin{ejem}
Consider the identity $\mathrm{Id}\colon  (\EE, c_1) \to (\EE, \bar c_1)$ between a
ring spectrum with two different orientations $c_1$ and $\bar c_1$. The explicit
computations of $G^{-1}$ is a classic subject on formal group laws (\emph{cf.}
\cite[\S 5.2]{Deglise} for a review in the context of the Riemann-Roch theorem). The
simplest example being
$$
\bar c_1 (L)=-c_1(L^*).
$$
Recall from the proof of Proposition \ref{Clase de Thom} that the canonical short
exact sequence $0\to \O_{\PP^1}(-1)\to \O_{\PP^1}\oplus \O_{\PP^1}\to Q \to 0$
satisfies that $Q=\O_{\PP^1}(1)$ so we have $c_1(\O_{\PP^1}(1)) =
-c_1(\O_{\PP^1}(-1))$ for any orientation. Therefore the class $\bar c_1
(\O_{\PP^\infty}(-1))$ as defined above is always an orientation. If $F$ is the
formal group law of $c_1$ then the series $G$ of Proposition \ref{serie c1} in this
case is the \emph{formal inverse} $\mu$ of $F$, i.e., the series satisfying
$F(x,\mu(x))=0$. However, it is much easier to compute Chern classes explicitly by
the splitting principle and the projective bundle theorem obtaining
$$
\bar c_i (E)=(-1)^i c_i(E^*) \qquad \mbox{and} \qquad  \sum_{i=0}^n(-1)^i
c_i(E^*)y^{n-i}=0 \in\EE(\PP(E))
$$
where $y=c_1(\O_{\PP(\EE)} (1))$.
\end{ejem}

The Riemann-Roch theorem for modules is a direct consequence of the case of rings.

\begin{teo}\label{Teo RR modulos}
Let $\varphi\colon  (\EE , c_1) \to (\FF , \bar c_1)$ be a morphism of oriented
absolute spectra. Denote $G\in \FF[S][[t]]$ the series such that $\varphi
(c_{1})=G(\bar c_1)\cdot \bar c_1\in \FF (\PP_S ^\infty)$ and $G^{-1}_\times$ be the
multiplicative extension of $G^{-1}$. Let $\Phi\colon \MM \to \MM'$ be a
$\varphi$-morphism of absolute modules and $f\colon Y\to  X$ be a projective lci
where $f=p\circ i$ and denote $T_f= i^{*}T_p-N_i\in K_0(Y)$ the virtual tangent
bundle. Then the diagram
$$
\xymatrix{
\MM(Y)\ar[r]^{f_*} \ar[d]_{G^{-1}_\times (T_f)\Phi_Y}& \MM(X)\ar[d]^{\Phi_X}\\
\MM' (Y) \ar[r]^{\bar f_*} & \MM' (X) }\phantom{G^{-1}_\times}
$$
commutes. In other words, for $m \in \MM(Y)$ we have
$$
\Phi_X(f_*(m))=\bar f_*\bigl(G^{-1}_\times(T_f)\cdot\Phi_Y(m)\bigr).
$$
\end{teo}
\qed

Recall that in Section \ref{Sec Ejemplos Modulos} we constructed different examples
of modules. Let us concretely remark the case of the relative cohomology.

\goodbreak

\begin{teo}\label{Teo RR motivico relativo}
Let $\varphi\colon \EE\to \FF$ be a morphism of strict oriented absolute ring
spectra. Let $g\colon T \to X$ be a morphism of schemes, $f\colon Y\to X$ be a
projective lci morphism and denote $g_{\scriptscriptstyle Y}\colon T\times_X Y\to Y$
and $T_f\in K_0(Y)$
the virtual tangent bundle of $f$. Assume in addition either $g$ is proper or $f$ is
smooth, then the diagram
$$
\xymatrix{ \EE(g_{{\scriptscriptstyle Y}})\ar[r]^{f_*}
\ar[d]_{G^{-1}_\times(T_f)\varphi}&
\EE(g)\ar[d]^{\varphi}\\
\FF(g_{{\scriptscriptstyle Y}}) \ar[r]^{ \bar f_*} & \FF(g) }
$$
commutes. In other words, for $m \in \EE(g_{\scriptscriptstyle Y})$ we have
$$
\varphi(f_*(m))=  f_*\bigl(G^{-1}_\times(T_f)\cdot \varphi(m)\bigr).
$$
\end{teo}
\demo We have a commutative diagram
$$
\xymatrix{\hofib _\EE(g) \ar[r]\ar@{-->}[d]_{\Phi_X} &\EE_X\ar[r]\ar[d]_{\varphi_X}
&Rg_*g^*\EE_X\ar[d]\\
\hofib _\FF(g)\ar[r] &\FF_X\ar[r]&Rg_*g^*\FF_X}
$$
so that there exists a $\varphi_X$-morphism of modules  $\Phi_X\colon\hofib
_\EE(g)\to \hofib _\FF(g)$. We conclude by applying Theorem \ref{Teo RR modulos} and
Proposition \ref{Propo cambio base liso hofib}.

\qed

\begin{ejem}
We review some concrete examples of this formula. Let $\varphi \colon \EE\to \FF$ a
morphism of strict absolute ring spectra:
\begin{itemize}
\item
Let $i\colon Z\to X$ be a closed immersion and consider the cartesian square
$$
\xymatrix{ P\ar[r]^-{i'} \ar[d]_{\pi'}& B_ZX\ar[d]^{\pi} \\
Z\ar[r]^-{i} & X}
$$
where $B_ZX$ denotes the blow-up of $Z$ in $X$. Recall from Example \ref{Ejem
cohomologia relativa} that $\EE(\pi)\simeq \EE(\pi')=\EE(P)/\EE(Z)$. We deduce from
the Riemann-Roch theorem for modules that the square
$$
\xymatrix{ \EE(P)/\EE(Z)\ar[r]^-{i_*} \ar[d]_{G^{-1}_\times( -N_{Z/X})\varphi}&
\EE(P)/\EE(Z) \ar[d]^{\varphi} \\
\FF(P)/\FF(Z)\ar[r]^-{\bar i_*} & \FF(P)/\FF(Z)}
$$
commutes. Note that $i_*(m)=m\cdot \bigl((i \circ \pi')^*\eta_Z^X\bigr)=m\cdot
c_n(N_{Z/X}) $. Recall from Corollary \ref{Tod cn} that
$c_n(N_{Z/X})=G^{-1}_\times(-N_{Z/X})\cdot \bar c_n(N_{Z/X})$ so the formula also
follows from the Riemann-Roch theorem for rings.

\item
Let $S=\Spec (k)$ be a point. Denote $g\colon X\to S$ be $k$-scheme and $f\colon Y\to
S$ be a smooth projective $k$-scheme. We have that $\EE(g)=\widetilde{\EE}(X) =
\EE(X)/\EE(S)$ and that $\EE(g_{{\scriptscriptstyle Y}})=\EE(X\times Y)/\EE(Y)$. Then
the square
$$
\xymatrix{ \EE(X\times Y)/\EE(Y)\ar[r]^-{f_*} \ar[d]_{G^{-1}_\times (T_f)\varphi}&
\widetilde{\EE}(X)\ar[d]^{\varphi} \\
\FF(X\times Y)/\FF(Y)\ar[r]^-{\bar f_*} & \widetilde{\FF}(X) }
$$
commutes.

Now assume both $\EE$ and $\FF$ satisfy the K\"{u}nneth formula. Then $\EE(X\times
Y)/\EE(Y) = \widetilde{\EE}(X)\otimes \EE(Y)$ and
$$
f_*=1\otimes f_*\colon \widetilde{\EE}(X)\otimes \EE(Y)\to \widetilde{\EE}(X)\otimes
\EE(S)=\widetilde{\EE}(X)
$$
and $\bar f_*=1\otimes \bar f_*$. In this case this formula also follows from the
Riemann-Roch theorem for ring spectra.

\item
\emph{Residual Riemann-Roch:} Let $i\colon Z\to X$ be a closed immersion of open
complement $j\colon U\to X$. Recall that $\EE(j)=\EE_Z(X)$. Note that the morphism of
modules $\varphi\colon \hofib_\EE(j)\to\hofib _\FF(j)$ fits into a morphism of
distinguished triangles. We deduce that the square
$$
\xymatrix{ \EE(U)\ar[r]^-{\delta} \ar[d]_{\varphi}& \EE_Z(X)\ar[d]^{\varphi} \\
\FF(U)\ar[r]^-{\bar \delta} & \FF_Z(X),}
$$
where $\delta$ denotes the connecting, commutes. If both $Z$ and $X$ are smooth then
we have the purity isomorphism $\EE(Z)\simeq \EE_Z(X)$ and we deduce D\'{e}glise's
residual Riemann-Roch formula (\cf \cite[4.2.3]{Deglise}): The square
$$
^{\phantom{G^{-1}_\times (_{Z/X})}}\xymatrix{ \EE(U)\ar[r]^-{\delta \p_i}
\ar[d]_{\varphi}& \EE(Z)\ar[d]^{G^{-1}_\times (-N_{Z/X}) \varphi} \\
\FF(U)\ar[r]^-{\bar \delta \bar \p_i} & \FF(Z)}
$$
commutes.

\item
Let $i\colon Z\to X$ be a closed immersion of open complement $j\colon U\to X$. We
have that $\EE(i)=\EE_c(U)$. Let $f\colon Y\to X$ be a projective lci morphism, the
square
$$
\xymatrix{ \EE_c(U\times_X Y)\ar[r]^-{f_*} \ar[d]_{G^{-1}_\times(T_f)\varphi}&
\EE_c(U)\ar[d]^{\varphi}
\\
\FF_c(U\times_X Y)\ar[r]^-{\bar f_*} & \FF_c(U),}
$$
commutes.

\end{itemize}
\end{ejem}

\subsection{Applications}

The main application we are interested in is the Grothendieck-Riemann-Roch theorem
for higher $K$-theory and the Riemann-Roch theorem for the relative cohomology of a
morphism. We afterwards review some other Riemann-Roch type formulas and the
arithmetic Riemann-Roch theorem. Recall from Example \ref{Ejem espectros en anillo}
that the Chern character is a morphism of strict absolute ring spectra $\ch\colon
\KGL_\QQ \to \bigoplus_i\HB (i)[2i]$. Denote $\Td$ the multiplicative extension of
the Todd series $\frac{t}{1-e^{-t}}$ (\cf Paragraph \ref{Tod}).

\goodbreak

\begin{teo}[Riemann-Roch]\label{Riemann Roch}
Let $f\colon Y\to  X$ be a projective lci where $f=p\circ i$ and denote $T_f\coloneqq
i^{*} T_p - N_i\in K_0(Y)$ the virtual tangent bundle. Then the diagram
$$
\xymatrix{
KH(Y)_\QQ\ar[r]^{f_*} \ar[d]_{\Td(T_f)\ch}& KH(X)_\QQ\ar[d]^{\ch}\\
H_{\M}(Y,\QQ) \ar[r]^{ f_*} & H_{\M}(X,\QQ) }
$$
commutes. In other words,
$$
\ch(f_*(a))= f_*\bigl(\Td(T_f)\cdot \ch(a)\bigr).
$$
\end{teo}
\demo The result follows from Theorem \ref{Motivic Riemann Roch} applied to $\ch$.
Recall that $\ch (L)=e^{c_1^{\HB}(L)}$, that $c_1^{\KGL} (L) = 1 - L^*$  and that
$c_1^{\HB}$ is additive so, in particular, $c_1^{\HB} (\O_{\PP^n}(1))=-c_1^{\HB}
(\O_{\PP^n}(-1))$. Denote $x =c_1^{\HB} = c_1^{\HB}(\O_{\PP^\infty}(-1))$ and $y =
c_1^{\KGL}=c_1^{\KGL}(\O_{\PP^\infty}(-1))$, we have
$$
\ch(y)=1-e^{-x}= x\cdot\frac{1-e^{-x}} {x}.
$$
Therefore $G=\frac{1-e^{-t}}{t}$ and $G^{-1}_\times=\Td$ the multiplicative extension
of $\frac{t}{1-e^{-t}}$.

\qed

\medskip

\begin{obser}
Note that in the case over a base field of exponential characteristic due to
comparison result of \cite{CD3} this theorem applies to \emph{cdh}-motivic
cohomology.
\end{obser}

A general Riemann-Roch statement as in \cite{Gillet} follows from the fact that
Beilinson motivic cohomology spectrum is universal for spectra with additive
orientations (\emph{cf.} \cite[14.2.16]{CD} and \cite[5.3.9]{Deglise}).

\begin{propo}\label{HB Universal}
Let $(\EE, c_1)$ be an oriented absolute ring spectrum in $\SH(S)_\QQ$ with $c_1$
additive. Then there exists a unique morphism of absolute spectra
$$
\varphi \colon  \HB \to \EE.
$$
Moreover, the morphism satisfies that $\varphi_{\PP^\infty} (c^{\HB}_1)=c_1 \in
\EE^{2,1}(\PP^\infty)$.
\end{propo}

Recall that some examples of oriented absolute ring spectra with additive
orientations are those coming from real absolute Hodge and Deligne-Beilinson
cohomology, rigid syntomic cohomology, and mixed Weil theories. Let now $S=\Spec (k)$
for $k$ a perfect field for mixed Weil theories, a field of characteristic zero for
real absolute Hodge and Deligne-Beilinson cohomology, or a residue field of a
$p$-adic field for rigid syntomic cohomology. The next result follows from the
Riemann-Roch theorem and Proposition \ref{HB Universal}:

\goodbreak

\begin{coro}
Let $\mathrm{H}$ denote either real absolute Hodge cohomology, real Deligne-Beilinson
cohomology, rigid syntomic cohomology or any cohomology coming from a mixed Weil
theory. Let $S$ be as above so that $\mathrm{H}$ is defined and let $f\colon Y\to X$
be a projective lci morphism of $S$-schemes. Then, with previous notations, the
diagram
$$
\xymatrix{
KH(Y)_\QQ\ar[r]^{f_*} \ar[d]_{\Td(T_f)\ch}& KH(X)_\QQ\ar[d]^{\ch}\\
\mathrm{H} (Y) \ar[r]^{ f_*} & \mathrm{H} (X) }
$$
commutes. In other words, for $a \in KH(Y)_\QQ$ we have
$$
\ch(f_*(a))= f_*\bigl(\Td(T_f)\cdot \ch(a)\bigr).
$$
\end{coro}
\qed

\medskip

Another general type of morphism of oriented absolute ring spectra to which the
motivic Riemann-Roch theorem applies are those coming from algebraic cobordism
$\MGL$. Recall that $\MGL$ is the universal oriented absolute ring spectrum (see
\cite{Vezzosi}).

\begin{propo}
Let $(\EE, c_1)$ be an oriented absolute ring spectrum. Then there exists a unique
morphism of absolute ring spectra
$$
\varphi\colon  \MGL \to \EE
$$
such that $\varphi (c_1 ^{\MGL})=(c_1) \in \EE ^{2,1}(\PP ^\infty)$.
\end{propo}

\qed

\medskip

Since this morphism preserves the orientation then Theorem \ref{RR Facil} applies to
them.

\begin{coro}
Let $(\EE, c_1)$ be an oriented absolute ring spectrum and $f\colon Y\to X$ be a
projective lci morphism. Then, with previous notations, for $a \in  \MGL(Y)$ we have
$$
\varphi_X(f_*(a))= f_*(\varphi_Y(a)).
$$
\end{coro}
\qed

\medskip

We apply the Riemann-Roch theorem for modules \ref{Teo RR modulos} to the examples we
described in Section \ref{Sec Ejemplos Modulos}.

\goodbreak

\begin{teo}
Let $g\colon T \to X$ be a morphism of schemes, $f\colon Y\to X$ be a projective lci
morphism and denote $g_{\scriptscriptstyle Y}\colon T\times_X Y\to Y$ and
$T_f\in K_0(Y)$ the virtual
tangent bundle of $f$. Assume in addition either $g$ is proper or $f$ is smooth, then
the diagram
$$
\xymatrix{
KH(g_{{\scriptscriptstyle Y}})_\QQ\ar[r]^{f_*} \ar[d]_{\Td(T_f)\ch}&
KH(g)_\QQ\ar[d]^{\ch}\\
H_{\M}(g_{{\scriptscriptstyle Y}},\QQ) \ar[r]^{f_*} & H_{\M}(g,\QQ) }
$$
commutes. In other words, for $m \in KH(g_{\scriptscriptstyle Y})_\QQ$ we have
$$
\ch(f_*(m))=  f_*\bigl(\Td(T_f)\cdot \ch(m)\bigr).
$$
\end{teo}

\qed

\medskip

We also obtain the arithmetic Riemann-Roch theorem of \cite{HS} as a consequence of
the Riemann-Roch theorem for modules.

\begin{teo}[Arithmetic Riemann-Roch]
Let $f\colon Y \to X$ be a projective morphism between smooth schemes over an
arithmetic ring and $T_f\in K_0(Y)$ the virtual tangent bundle. Then the diagram
$$
\xymatrix{ \hatKH(Y)_\QQ\ar[r]^{f_*} \ar[d]_{\Td(T_f)\hatch}&
\hatKH(X)_\QQ  \ar[d]^{\hatch}\\
\widehat {H}_{\M}(Y,\QQ) \ar[r]^{ f_*} & \widehat{H}_{\M}(f,\QQ) }
$$
commutes. In other words, for $m \in \hatKH(Y)_\QQ$ we have
$$
\hatch(f_*(m))=  f_*\bigl(\Td(T_f)\cdot \hatch(m)\bigr).
$$
\end{teo}

\qed

\appendix

\section{Absolute Hodge cohomology} \label{Hodge absoluto}

In this Appendix we apply a theorem of D\'{e}glise and Mazzari to give direct
construction of the (real) absolute Hodge spectrum representing absolute Hodge
cohomology with real coefficients with no smoothness assumption. Since in \cite{HS}
the authors asked explicitly if the Deligne-Beilinson spectrum represented the
Deligne-Beilinson cohomology on singular schemes we also prove it for the
Deligne-Beilinson spectrum.

We refer to Drew's thesis \cite{Drew} for the original construction of the (rational)
absolute Hodge spectrum and a more complete treatment on the subject. We refer to
\cite{Beilinson}, \cite{Jannsen}, \cite{Burgos3} and \cite{Beilinson2}, \cite{EV},
\cite{Burgos2} for more details of the following constructions for absolute Hodge and
Deligne-Beilinson cohomology respectively.

\begin{parra}
Let $X$ be a smooth complex variety. We can find a proper complex variety $\bar X$
and an open embedding $j\colon X \to \bar X$ such that $D=\bar X - X$ is a normal
crossing divisor. We denote by $A_{\bar X}^* (\log D)$ the complex of smooth
differential forms with logarithmic singularities along $D$ (\emph{cf.}
\cite{Burgos}). Taking limit over all suitable compactifications we define
$$
A_{\log } ^*(X)=\lim A_{\bar X}^* (\log D)
$$
the complex of smooth differential forms with logarithmic singularities along
infinity. The complex $A_{\log } ^*(X)$ has a natural filtration $W$ which assigns
weight zero to the sections of $A^*(X)$ and weight one to the sections $\d z_i/z_i$
and $\d \bar z_i/z_i$. The complex $A_{\log } ^*(X)$ also has a natural Hodge
filtration $F$, as well as a subcomplex $A_{\log ,\RR}^*(X)$ of differential forms
invariant under complex conjugation. Therefore it defines an $\RR$-Hodge complex. The
absolute Hodge cohomology of $X$ with real coefficients is defined as
$$
H_{\mathrm{AH}}^p(X, \RR (q))= H^p (\tilde \Gamma (A_{\log}(X),q))
$$
where
$$
\tilde \Gamma (A_{\log}(X),q)= \mathrm{cone}( (2\pi i)^q \hat W_{2q} A_{\log ,
\RR}(X) \oplus \hat W_{2q}\cap F^q A_{\log}(X) \to \hat W_{2q}A _{\log} (X))[-1]
$$
and $\hat W$ denotes the decal\'{e} filtration of $W$ (\emph{cf.} \cite[1.1.2]{Deligne}).

The \emph{real Deligne-Beilinson cohomology} is obtained by ignoring the weight
filtration. That is to say, we define it as
$$
H_{\mathrm{D}}^p(X,\RR (q))=H^p(\Gamma (A_{\log}(X) ,q))
$$
where
$$
\Gamma(A_{\log}(X),q)=\mathrm{cone}( (2\pi i)^q A_{\log , \RR}(X) \oplus F^q
A_{\log}(X) \to A _{\log} (X))[-1].
$$

Both the real absolute Hodge and the Deligne-Beilinson cohomology can also be
computed by means of the \emph{Thom-Whitney simple} introduced in \cite{Navarro
Aznar}. Following \cite{Burgos3}, the Thom-Whitney simple has a concrete description
in these cases. Denote $L_1^*$ the differential graded commutative $\RR$-algebra of
algebraic forms over $\AA^1_\RR$, then the Thom-Whitney simple $\tilde \Gamma
_{\mathrm{TW}}(A_{\log}(X),q)$ for the real absolute Hodge cohomology is the
subcomplex of
$$
(2\pi i)^q \hat W_{2q} A_{\log , \RR}(X) \oplus \hat W_{2q}\cap F^q A_{\log}(X)
\oplus (L_1^*\otimes \hat W_{2q}A _{\log} (X))
$$
made by elements $(a,b,\omega)$ such that $\omega (0)=a$ and $\omega (1)=b$. The
differential is the natural one on each summand. These complex satisfy that
$$
H_{\mathrm{AH}}^p(X, \RR (q))=H^p(\tilde \Gamma _{\mathrm{TW}}(A_{\log}(X),q)).
$$
\end{parra}

\begin{defi}\label{Defi cuerpo aritmetico}
An \textbf{arithmetic field} is a triple $(k , \, \Sigma , \, \Fr )$ where $k$ is a
field, $\Sigma=\{ \sigma_1 , \ldots  , \sigma_n\}$ is a set of embeddings $k\to \CC$
and $\Fr \colon \CC^\Sigma \to \CC ^\Sigma$ is a $\CC$-antilinear involution such
that the image of $k$ in $\CC^\Sigma$ by $\sigma_1 \times \cdots \times\sigma_n$ is
invariant under $\Fr $. We call \emph{Frobenius} to the map $\Fr$. If $X$ is a
$k$-scheme we write $ X_\Sigma=\coprod \big ( X\times_{\sigma_i} \, \Spec k\big )$,
which is naturally a complex variety.
\end{defi}

Let $X$ be a smooth scheme over an arithmetic field. The Frobenius $\Fr $ induces a
$\CC$-linear action on $A_{\log }(X_\Sigma)$ by taking the action it induces on a
compactification $\bar X_\Sigma$ together with the complex conjugation. That is to
say, $\Fr f(x)=\bar f (\Fr (x))$. This action is compatible with the weight and Hodge
filtration. We consider
$$
\tilde \Gamma (A_{\log} (X_\RR),q)=\tilde \Gamma (A_{\log} (X_\Sigma),q)^{\Fr } \ \ ,
\ \  \Gamma (A_{\log }(X_\RR), q)=\Gamma (A_{\log }(X_\Sigma),q)^{\Fr }.
$$
We denote the cohomology they define as
$$
H_{\mathrm{D}}^p(X_\RR,\RR (q))=H^p(\tilde \Gamma (A_{\log} (X_\RR),q))\ \ , \ \
H_{\mathrm{AH}}^p(X_\RR, \RR (q))= H^p(\Gamma (A_{\log }(X_\RR), q)).
$$
As before, their respective Thom-Whitney simple also compute the cohomology. In the
case of the complex for absolute Hodge cohomology we denote it $\tilde \Gamma
_{\mathrm{TW}} (A_{\log} (X_\RR),q)$.

In \cite{HS} Holmstrom and Scholbach proved that there exist an absolute spectrum
$\EE_{\mathrm{D}} \in \SH(S)_\QQ$ that represents the Deligne-Beilinson cohomology
for smooth schemes. In other words, there exist an absolute spectrum
$\EE_{\mathrm{D}}$ such that for every $X$ smooth and every integers $p,q$ we have
$$
\EE_{\mathrm{D}}^{p,q}(X)=H_{D}^p(X,\RR(q)).
$$
This argument has been vastly generalized by D\'{e}glise and Mazzari in \cite[1.4.10]{DM}
by giving sufficient conditions for a family of presheaves $X\mapsto \F_i(X)$ for
$i\in \NN$ so that the cohomology they define is represented by an absolute ring
spectrum. That is to say, they give sufficient conditions on $(\F_i)_{i\in\NN}$ so
that there exist an absolute ring spectrum $\EE_\F$ satisfying
$$
H^p(\F_{q}(X))=\EE_\F ^{p,q}(X).
$$

\begin{parra}\label{Brown-Gersten}
In order to check that the family $X\mapsto \tilde \Gamma_{TW} (A_{\log}^*(X_\RR),i)$
satisfy the hypothesis of \emph{loc. cit.} let us introduce some notation. Consider
$\CC$ as an arithmetic field with $\Sigma=\{\mathrm{Id}, \, \sigma \}$ where $\sigma$
denotes the complex conjugation. Denote $c\colon \RR \to \tilde \Gamma_{TW}
(A_{\log}^*(\GG_{m \, \RR}),1)[1]$ the section given by
$$
\big (\frac{\d z}{z}+ \frac{\d \bar z}{\bar z},\, (\frac{\d z}{z}-\frac{\d \bar
z}{\bar z})i,\, (1-x)(\frac{\d z}{z}+\frac{\d \bar z}{\bar z})+x (\frac{\d
z}{z}-\frac{d \bar z}{\bar z})i+(\ln z\bar z +i\ln \frac{z}{\bar z})\d x \big )
$$
where the first term is in $ (2\pi i) \hat W_{2} A_{\log 0, \infty ,
\RR}((\PP^1)_\RR)$, the second belongs to $ \hat W_{2}\cap F^1 A_{\log 0 ,
\infty}((\PP^1)_\RR)$ and the third belongs to $L_1^*\otimes \hat W_{2}A _{\log 0 ,
\infty} ((\PP^1)_\RR)$. For a general arithmetic field we still denote $c$ the
section defined by taking $c$ on each component of $\GG_{m\, \Sigma}$.

Also recall that a \emph{distinguished square} is a commutative cartesian diagram
$$
\xymatrix{Y' \ar[r] \ar[d]& X' \ar[d]\\
Y \ar[r] & X}
$$
in $\Sm/S$ such that $Y\to X$ is an open immersion and $X'\to X$ is \'{e}tale and induces
an isomorphism $(X'\setminus Y')_{\mathrm{red}}\to (X\setminus Y)_{\mathrm{red}}$. We
say that a complex of presheaves of $R$-modules $\F$ on $\Sm/S$ has the
\emph{Brown-Gersten property with respect to the Nisnevich topology} if for every
distinguished square the diagram
$$
\xymatrix{\F(X) \ar[r] \ar[d]& \F (X') \ar[d]\\
\F(Y) \ar[r] & \F (Y')}
$$
is a homotopy pullback square in the category of complexes of $R$-modules.
\end{parra}

\begin{propo}
Let $k$ be an arithmetic field and denote $S=\Spec (k)$. Consider the family of
presheaves $X\mapsto \tilde \Gamma_{TW} (A_{\log}^*(X_\RR),i)$ on $\Sm/S$ together
with the section $c\colon \RR \to \tilde \Gamma_{TW} (A_{\log}^*(\GG_{m\,
\RR}),1)[1]$ defined above:

\begin{enumerate}

\item
They form an $\NN$-commutative graded monoid (\emph{cf.} \cite[1.4.9]{DM}) in the
category of complexes of $\RR$-linear presheaves on the category of affine and smooth
$S$-schemes.

\item
They have the Brown-Gersten property with respect to the Nisnevich topology.

\item
They are homotopy invariant, i.e.,
$H_{\mathrm{AH}}^p(\AA^1_X,\RR(q))=H_{\mathrm{AH}}^p(X,\RR(q))$.

\item
If $\bar c\in H_{\mathrm{AH}}^1(\GG_{m\, \RR},\RR(1))$ denotes the class of $c$, then
for any smooth scheme $X$ and any integers $p$, $q$ the map
$$
\begin{matrix}
H_{\mathrm{AH}}^p(X_\RR,\RR(q))&\longrightarrow  & H_{\mathrm{AH}}^{p+1}((X\times
\GG_m)_\RR,\RR(q+1))/H_{\mathrm{AH}}^p(X_\RR,\RR(q)) \\
 x &\mapsto & [x\times \bar c]
\end{matrix}
$$
where $[\raya ]$ denotes the class defined in the quotient, is an isomorphism.

\item
If $u\colon \GG_m \to \GG_m$ is the inverse map of the group scheme $\GG_m$ and $\bar
c'$ is the image of $c$ in $H^1_{\mathrm{AH}}(\GG_{m\, \RR},\RR(1))$ then $u^*(\bar
c')=-\bar c'$.

\end{enumerate}
\end{propo}
\demo Although these properties are well known for experts let us review them for the
sake of completeness. The Thom-Whitney simple has a well defined associative and
commutative product (\emph{cf.} \cite[\S 6]{Burgos}), from which point 1 follows.

For point 2, first notice that the Brown-Gersten property is stable by direct sums
and cones of maps. Therefore it is enough to prove it for $ \hat W_{q} A_{\log ,
\RR}$, $ \hat W_{2q} \cap F^q A_{\log}$ and $\hat W_{q}A _{\log}$. The \'{e}tale descent
for $A_{\log , \RR} $, $ F^q A_{\log}$ and $A _{\log} $ may be found in
\cite[2.9]{HS}, from which the Brown-Gersten property follows. For the weight
filtration, consider a distinguished square as in \ref{Brown-Gersten}. We have a
short exact sequence
$$
0\to A_{\log}(X)\to A_{\log}(X')\bigoplus A_{\log}(Y)\to A_{\log}(Y') \to 0
$$
as well as for $A_{\log, \RR}$ and $ F^q A_{\log}$. These morphisms are strict both
for the Hodge and the decal\'{e} Weight filtration, so the Brown-Gersten property readily
follows for $ \hat W_{q} A_{\log , \RR}$, $ \hat W_{2q} \cap F^q A_{\log}$ and $\hat
W_{q}A _{\log}$. We conclude by taking invariants on the action induced by the
Frobenius.

For point 4 note that the absolute Hodge cohomology of $\GG_m$ is null apart from the
groups $H^0(\GG_{m},\RR(0))=\RR$ and $H^1(\GG_{m},\RR(1))=\RR$. It is easy to see
from the definition of absolute Hodge cohomology that group $
H_{\mathrm{AH}}^{p+1}((X\times \GG_m)_\RR,\RR(q+1))$ equals
$$
\big (H_{\mathrm{AH}}^{p} (X_\RR,\RR(q))\otimes
H_{\mathrm{AH}}^{1}(\GG_{m\,\RR},\RR(1)) \big )\oplus \big ( H_{\mathrm{AH}}^{p+1}
(X_\RR,\RR(q+1))\otimes H_{\mathrm{AH}}^{0}(\GG_{m\,\RR},\RR(0))\big )
$$
from which point 4 follows. Finally, a direct computation concludes point 5.

\qed

\begin{coro}
Denote $\EE_{\mathrm{AH}}$ the oriented absolute ring spectrum  constructed in
\cite[1.4.10]{DM} out of the presheaves $(\tilde \Gamma_{TW}(A_{\log}^*(\raya
_\RR),i))_{i\in \NN}$, which we call the absolute Hodge spectrum. Then
$\EE_{\mathrm{AH}}$ represents real absolute Hodge cohomology on smooth schemes. In
other words, for any smooth $S$-scheme $X$ and any $p,q\geq 0$
$$
\EE^{p,q}_{\mathrm{AH}}=H^p_{\mathrm{AH}}(X_\RR,\RR(q)).
$$
\end{coro}
\qed

Let us now prove that the absolute Hodge spectrum represents absolute Hodge
cohomology for general schemes. The same method will apply also for the
Deligne-Beilinson cohomology.

\begin{parra}
Let $Z$ be complex variety, following \cite[8.3]{Deligne2} we can find a diagram
$\bar X _\bullet \hookleftarrow X_\bullet \buildrel{p}\over{\to} Z$ so that
$X_\bullet$ and $\bar X_\bullet$ are simplicial complex varieties,  $p$ satisfies
cohomological descent (in particular, it is a hypercovering for the
\emph{h}-topology), $\bar X_\bullet$ is proper smooth and $\bar X_\bullet -
X_\bullet$ is a normal crossing divisor. If $\bar X'_\bullet \hookleftarrow
X'_\bullet\to Z$ is a second diagram we have the isomorphisms
$H_{\mathrm{AH}}(X_\bullet,\RR)\simeq H_{\mathrm{AH}} (X'_\bullet ,\RR)$ and
$H_{\mathrm{D}}(X_\bullet,\RR)\simeq H_{\mathrm{D}}(X'_\bullet,\RR)$ for both the
absolute Hodge and Deligne-Beilinson cohomology of those simplicial varieties.
Therefore, we call the absolute Hodge and the Deligne-Beilinson cohomology of $Z$ to
$$
H_{\mathrm{AH}}(Z,\RR)\coloneqq H_{\mathrm{AH}}(X_\bullet,\RR) \qquad
H_{\mathrm{D}}(Z,\RR)\coloneqq H_{\mathrm{D}}(X_\bullet,\RR)
$$
This construction is compatible with the Frobenius action so we define the groups
$H_{\mathrm{AH}}(Z_\RR,\RR)\coloneqq H_{\mathrm{AH}}(X_{\bullet \, \RR},\RR)$ and
$H_{\mathrm{D}}(Z_\RR,\RR)\coloneqq H_{\mathrm{D}}(X_{\bullet \, \RR},\RR)$.

Recall that any rational oriented absolute ring spectrum is also a Beilinson motive
(\emph{cf.} \cite[14.2.16]{CD} ). The category of Beilinson motives $\DMb (S)$
satisfies the \emph{h}-descent (\emph{cf.} \cite[3.1]{CD}). This implies that if $X$
is a scheme and $X_\bullet\to X$ is a \emph{h}-hypercover then for any oriented
absolute ring spectrum we have
$$
\EE (X_\bullet) = \EE (X),
$$
where $\EE (X_\bullet)$ denotes the cohomology of the simplicial scheme (\emph{cf.}
\cite[\S 3.1]{CD} or \cite[\S 2.2.1]{DM} for $\DMb$).
\end{parra}

\begin{coro}
Let $\EE_{\mathrm{AH}}$ and $\EE_{D}$ be the real absolute Hodge and
Deligne-Beilinson spectra, then for any scheme $Z$ over an arithmetic field and every
$p$, $q\geq 0$ we have
$$
\EE_{\mathrm{AH}}^{p,q}(Z)=H_{\mathrm{AH}}^p(Z_\RR,\RR (q))\quad \mbox{and} \quad
\EE_{\mathrm{D}}^{p,q}(Z)=H_{\mathrm{D}} ^p(Z_\RR,\RR (q)).
$$
\end{coro}

\end{document}